\documentclass[twocolumn]{autart} 

\usepackage{graphicx} 
\usepackage{enumerate}
\usepackage[round, semicolon]{natbib}
\usepackage[letterpaper, left=0.75in, right=0.75in, top=1in, bottom=1in]{geometry}
\usepackage{subcaption}
\usepackage{amsmath} 
\usepackage{amssymb}  
\usepackage{color}
\usepackage[T1]{fontenc}
\usepackage{todonotes}
\usepackage{comment}
\usepackage{wrapfig}

\newtheorem{innercustomthm}{Theorem}
\newenvironment{customthm}[1]
  {\renewcommand\theinnercustomthm{#1}\innercustomthm}
  {\endinnercustomthm}

\newtheorem{innercustomlem}{Lemma}
\newenvironment{customlem}[1]
  {\renewcommand\theinnercustomlem{#1}\innercustomlem}
  {\endinnercustomlem}

\newtheorem{innercustomass}{Assumption}
\newenvironment{customass}[1]
  {\renewcommand\theinnercustomass{#1}\innercustomass}
  {\endinnercustomass}  
  
\newtheorem{innercustomcor}{Corollary}
\newenvironment{customcor}[1]
  {\renewcommand\theinnercustomcor{#1}\innercustomcor}
  {\endinnercustomcor}  
  
\newtheorem{innercustomrem}{Remark}
\newenvironment{customrem}[1]
  {\renewcommand\theinnercustomrem{#1}\innercustomrem}
  {\endinnercustomcor}   
 
\newtheorem{innercustomcase}{Case}
\newenvironment{customcase}[1]
  {\renewcommand\theinnercustomcase{#1}\innercustomcase}
  {\endinnercustomcase}

\begin{document}

\begin{frontmatter}

\title{
Conditions to Provable System-Wide Optimal Coordination of Connected and Automated Vehicles\thanksref{footnoteinfo}}

\author[UD]{A M Ishtiaque Mahbub}\ead{mahbub@udel.edu},
\author[UD]{Andreas A. Malikopoulos}\ead{andreas@udel.edu}

\address[UD]{Department of Mechanical Engineering, University of Delaware, 126 Spencer Lab, 130 Academy St, Newark DE 19716}

\thanks[footnoteinfo]{This research was supported by ARPAE's NEXTCAR program under the award number DE-AR0000796.}

%
\begin{abstract}
Connected and automated vehicles (CAVs) provide the most intriguing opportunity to improve energy efficiency, traffic flow, and safety. In earlier work, we addressed the constrained optimal coordination problem of CAVs at different traffic scenarios using Hamiltonian analysis. In this paper, we investigate the properties of the unconstrained problem and provide conditions under which different combination of the state and control constraints become active. We present a condition-based computational framework that improves on the standard iterative solution procedure of the constrained Hamiltonian analysis. Finally, we derive a closed-form analytical solution of the constrained optimal control problem and validate the proposed framework using numerical simulation. The solution can be derived without any recursive steps, and thus it is appropriate for real-time implementation on-board the CAVs.
\end{abstract}
\end{frontmatter}

\section{Introduction} \label{sec:1}
\subsection{Motivation} The implementation of an emerging transportation system with connected and automated vehicles (CAVs) enables a novel computational framework to provide real-time control actions that optimize energy consumption and associated benefits. From a control point of view, CAVs can alleviate congestion at different traffic scenarios, reduce emission, improve fuel efficiency and increase passenger safety; see \cite{Margiotta2011, Malikopoulos2017}. Urban intersections, merging roadways, highway on-ramps, roundabouts and speed reduction zones along with the driver responses to various disturbances are the primary sources of bottlenecks that contribute to traffic congestion; see \cite{Malikopoulos2013}.  

\subsection{Literature Review}
Several research efforts have used optimal control theory to investigate how CAVs can potentially improve energy efficiency and travel time in these traffic scenarios.
Early efforts reported in \cite{Levine1966} and \cite{Athans1969} considered a single string of vehicles that was coordinated through a traffic conflict zone with a linear optimal regulator. \cite{Shladover1991} discussed the lateral and longitudinal control of CAVs for the automated platoon formation. \cite{Varaiya1993} outlined the key features of an automated intelligent vehicle/highway system, and proposed a basic control system architecture. \cite{Dresner2004} proposed the use of the reservation scheme to control a signal-free intersection of two roads. Since then, several research efforts have considered reservation approaches for coordination of CAVs at urban intersections; see \cite{Dresner2008, DeLaFortelle2010, Huang2012, Au2010a}. \cite{Alonso2011} proposed a control framework where a CAV can derive its safe crossing schedule to avoid collision with a human-driven vehicle.
Several approaches for coordinating CAVs that have been reported in the literature have proposed the use of centralized control, {where there is
at least one task in the system that is globally decided for all
vehicles by a single central controller}; see \cite{Dresner2008, DeLaFortelle2010, Huang2012, lu2003longitudinal, xu2018cooperative, bakibillah2019optimal}. {Some approaches have focused on coordinating CAVs at intersections to improve traffic flow; see \cite{Yan2009,kim2014}, or travel time; see \cite{Raravi2007}, while other approaches have focused on energy consumption improvement; see \cite{mahler2014optimal,sciarretta2015optimal,wan2016optimalmixed}.}

{
Some optimal control approaches reported in the literature have used standard Hamiltonian analysis for CAV control and coordination, e.g., \cite{Zhao2019CCTA-1, wang2019intersection}; while other approaches have employed model predictive control; see \cite{kim2014,Makarem2012}. Dynamic programming (DP) has also been used to compute the optimal control input for CAVs, e.g., \cite{ozatay2017velocity}, \cite{mahler2014optimal}, and \cite{pei2019cooperative}. DP, however, may not be feasible for real-time implementation due to its high required computational effort.}
In optimal control approaches, the problem formulation may have different objective functions including vehicle travel time, e.g., \cite{Raravi2007}, energy consumption, e.g., \cite{sciarretta2015optimal}, passenger comfort, e.g., \cite{Ntousakis:2016aa}, etc. \cite{Raravi2007} formulated an optimization problem the solution of which aims at finding the minimum time once the merging sequence is determined. \cite{Kamal2013} proposed numerical algorithms based on Pontryagin's minimum principle for CAV coordination in a signal-free intersection. A virtual platoon-based cooperative control approach was discussed in \cite{huang2019cooperative} for on-ramp coordination. A hierarchical control framework using an upper-level CAV coordination and a low-level multiobjective optimization scheme was proposed in \cite{qian2015decentralized}. A similar hierarchical control framework has been reported by \cite{bakibillah2019optimal}, where a two-level combinatorial optimization problem is formulated for a cloud-based roundabout coordination system.

{In optimal control approaches, one key challenge is to handle the associated state, control and safety constraints.
\cite{min2019constrained} considered a platoon-based approach to coordinate CAVs through a merging roadway, and solved the constrained optimization problem with distributed model predictive control.
\cite{sciarretta2015optimal} developed an eco-driving controller for CAVs for adaptive cruise control maneuver, where the optimal control problem minimizes the energy consumption with speed constraint. \cite{wan2016optimalmixed} proposed a speed advisory system to minimize fuel consumption without considering the state and control constraints. \cite{han2018safe} proposed a safety based eco-driving control for the CAVs.
\cite{wang2019intersection} formulated the multi-objective optimization problem for the CAVs approaching intersection, and derived the analytic solution based on the Pontraygin's minimum principle.
\cite{ozatay2017velocity} provided a speed profile optimization framework for minimizing fuel consumption without considering any safety or acceleration/deceleration constraints.}

Recently, a decentralized optimal control framework was presented for coordinating CAVs in real time at different traffic scenarios such as on-ramp merging roadways, roundabouts, speed reduction zones and signal-free intersections; see \cite{Malikopoulos2017,mahbub2020decentralized, Malikopoulos2018c,mahbub2020ACC-2}. This framework uses a hierarchical structure consisting of an upper-level vehicle coordination problem to minimize travel time, and a low-level optimal control problem to minimize the energy of individual CAVs. 
A complete, analytical solution of the low-level control problem that includes the rear-end safety constraint, where the safe distance is a function of speed, was discussed in \cite{malikopoulos2019ACC, Malikopoulos2020}. A problem formulation for the upper-level optimization in which there is no duality gap, implying that the optimal time trajectory for each CAV does not activate any of the state, control, and safety constraints of the low-level optimization was presented in \cite{Malikopoulos2019CDC, Malikopoulos2020}.

Detailed discussions of the research efforts reported in the literature to date on coordination of CAVs can be found in recent survey papers; see \cite{Malikopoulos2016a,Guanetti2018}.

\subsection{Objectives and Contributions of the Paper}
{The standard methodology to solve the low-level optimal control problem; see \cite{Malikopoulos2017}; is to employ Hamiltonian analysis with interior point state and/or control constraints. Namely, we first start with the unconstrained arc and derive the solution of the low-level optimal control problem. If the solution violates any of the state or control constraints, then the unconstrained arc is pieced together with the arc corresponding to the violated constraint. The two arcs yield a set of algebraic equations which are solved simultaneously using the boundary conditions and interior constraints between the arcs. If the resulting solution, which includes the determination of the optimal
switching time from one arc to the next one, violates another constraint, then the last two arcs are pieced together with the arc corresponding to the new violated constraint, and we re-solve the problem with the three arcs pieced together. The three arcs will yield a new set of algebraic equations that need to be solved simultaneously using the boundary conditions and interior constraints between the arcs. The resulting solution includes the optimal switching time from one arc to the next one. The process is repeated until the solution does not violate any other constraints. This recursive process of piecing the arcs together to derive the optimal solution of the low-level problem can be computationally expensive and might prevent real-time implementation.}

{In this paper, we provide an in-depth analysis of different state and control constraint activation cases, and establish a rigorous framework that yields a closed-form analytical solution for the low-level optimal control problem formulation without requiring the recursive process described above. Thus, the proposed framework is appropriate for real-time implementation on-board the CAVs; see \cite{mahbub2020sae-1}.
The objectives of this paper are
(i) to derive a priori the different state and control constraint activation cases through a rigorous mathematical analysis, 
(ii) to simplify the recursive process required to derive the optimal constrained solution of the Hamiltonian analysis for the low-level optimal control problem, and
(iii) to increase the computational efficiency of the derivation of the solution in (i) by eliminating numerical computations.}
    
Thus, the contributions of this paper are: 
(1) an in-depth exposition of the properties of the different combinations of the state and control constraint activation cases and a set of a priori conditions to identify the constrained solution without any recursive steps, and
(2) an explicit expression of the junction point between the constrained and unconstrained arcs leading to a closed-form analytical solution of the constrained optimal control problem.
%
{In earlier work, we reported a limited-scope analysis along with some preliminary results about the conditions for state and control constraint activation; see \cite{Mahbub2020ACC-1}.}

\subsection{Comparison With Related Work}
The framework that we report in this paper advances the state of the art in the following ways.
First, the solution to the state and control unconstrained control problem presented in \cite{Malikopoulos2018c} and \cite{Ntousakis:2016aa} shows acceleration spikes (jerk) at the boundaries of the optimization horizon, possibly exceeding the vehicle's physical limitation and giving rise to undesired driving experience. { In addition, the unconstrained solution can only guarantee that none of the constraints are violated at the boundaries of the optimization horizon only. 
In our proposed framework, we can guarantee that none of the the state and control constraints are violated throughout the entire optimization horizon.}
%
%
Second, in contrast to some approaches reported in the literature, e.g., \cite{wan2016optimalmixed}, \cite{ozatay2017velocity} and \cite{han2018safe}, where either the state or the control constrained optimal control problem was addressed, our framework addresses all state and control constraints cases. Moreover, we explicitly include the state and control constraints in the Hamiltonian analysis as opposed to using a feasibility zone; see \cite{wang2019intersection}.
%
Third, several approaches have considered free terminal time to address the state/control constraints within the optimization horizon; see \cite{wang2019intersection, zhang2019decentralized}.
In contrast, in our framework, we incorporate the constraints in the low-level control problem with the fixed time horizon.
%
Fourth, the solution of the constrained optimal control problem requires piecing the unconstrained and constrained arcs together resulting in recursive numerical computations until all of the constraint activation cases are resolved; see \cite{Malikopoulos2017}, \cite{malikopoulos2019ACC} and \cite{zhang2019decentralized}. In our proposed framework, we eliminate this recursive procedure to derive a real-time implementable closed-form analytical solution.
%
%
Finally, the solution of the constrained optimization problem using Hamiltonian analysis reported in some approaches, e.g., \cite{Malikopoulos2017}, \cite{malikopoulos2019ACC} and \cite{zhang2019decentralized}, only addresses different constraint activation cases without addressing the explicit interdependence between multiple constraint activation. In this paper, we explore the interdependence of the combination of the constraint activation cases and explicitly provide the conditions for their realization.
\subsection{Organization of the paper}
The remainder of the paper is organized as follows. In Section II, we introduce the problem formulation and present the unconstrained case. In Section III, we discuss different aspects of the state and control constrained formulation in detail. In Section IV, we provide the closed-form analytical solution of the constrained optimal control problem. In Section V, we evaluate the effectiveness of the proposed approach in a simulation environment. Finally, we draw concluding remarks and discuss potential directions for future research in Section VI. 

\section{Problem Formulation} \label{sec:2}
\begin{figure}[ht]
\centering
\includegraphics[width=3.3 in]{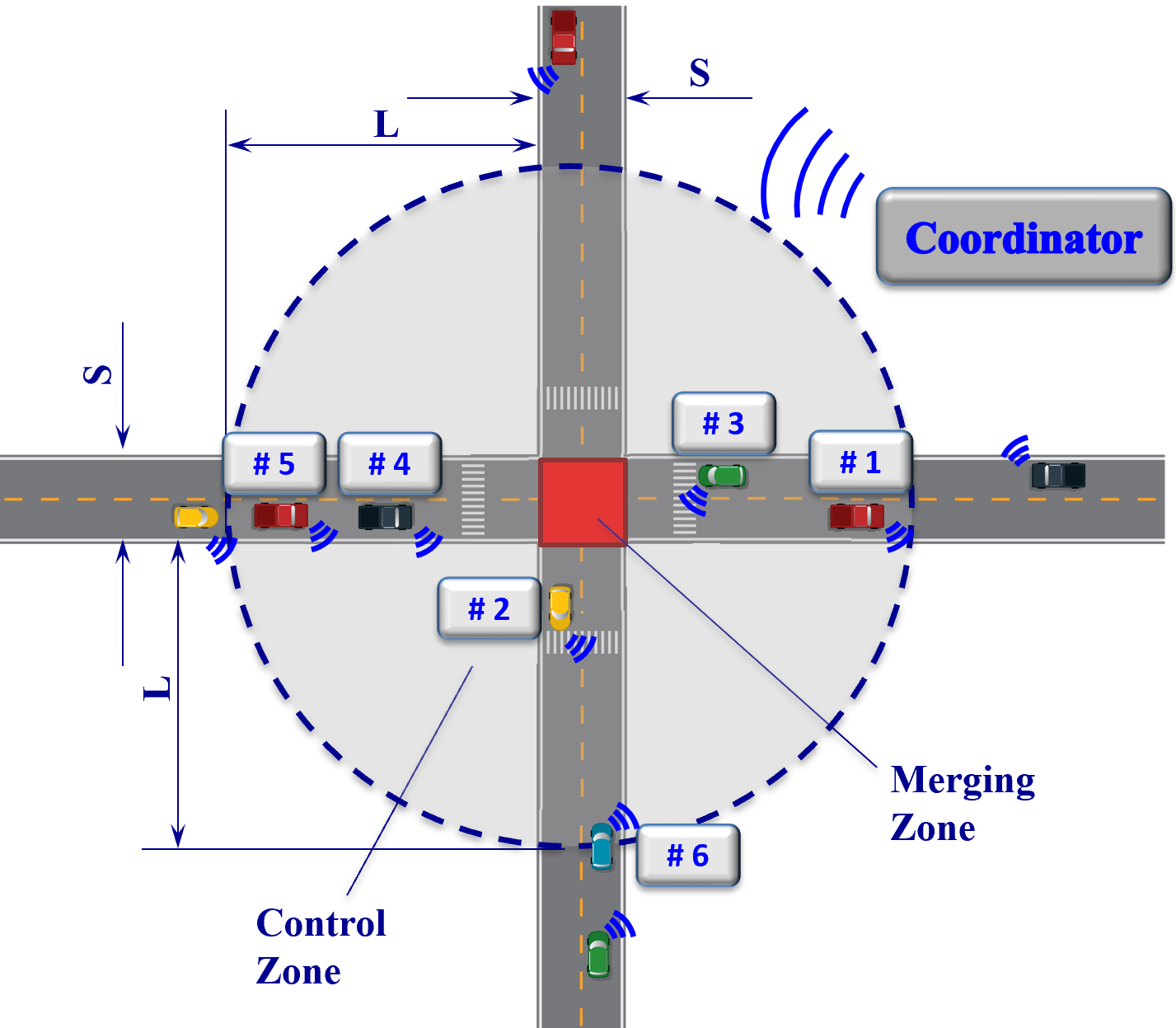} 
\caption{A traffic network of connected automated vehicles approaching a four-way signal-free intersection.}%
\label{fig:intersection}%
\end{figure}

We consider CAVs {travelling} through a traffic network containing a four-way signal-free intersection, as shown in Fig. \ref{fig:intersection}. Although our analysis can be applied to any traffic scenario, e.g., {merging at roadways, roundabouts, and passing through speed reduction zones,} we use an
intersection (Fig. \ref{fig:intersection}) as a reference to present the fundamental ideas and results of this paper, since an intersection provides unique features
making it technically more challenging compared to
other traffic scenarios. We define the area illustrated by the red square of dimension $S$ in Fig. \ref{fig:intersection} as the \emph{merging zone} where potential lateral collision of CAVs may occur. {Upstream of the merging zone, we define a \textit{control zone} of length $L$ inside of which CAVs can communicate with each other using a vehicle-to-vehicle communication protocol; see \cite{mahbub2020sae-1}. The intersection also has a \textit{coordinator} that communicates with the CAVs traveling inside the control zone.} Note that, the coordinator does not make any decisions for the CAVs. When a CAV enters the control zone, the coordinator receives its information and assigns a unique identity $i\in\mathbb{N}$ to it. 
{Let $\mathcal{N}(t)=\{1,\ldots, N(t)\}$, where $N(t)\in\mathbb{N}$ is the number of CAVs inside the control zone at time $t\in\mathbb{R}^{+}$, be the queue of CAVs to enter the merging zone shown in Fig. \ref{fig:intersection}. The time that a CAV $i\in \mathcal{N}(t)$ enters the control and merging zones is denoted by $t_i^0$ and $t_i^{m}$, respectively, while the time that a CAV $i$ exits the merging zone is denoted by $t_i^{f}$.}
{In our exposition, we assume that the queue $\mathcal{N}(t)$ and the optimal time to enter the merging zone $t_i^m$ is given a priori and can be derived by solving an upper-level vehicle coordination problem subject to rear-end and lateral safety constraints, as detailed in \cite{Malikopoulos2017, Mahbub2019ACC, mahbub2020decentralized}. Given $t_i^m$ a priori, the objective of each CAV $i\in\mathcal{N}(t)$ is to derive its optimal control input (acceleration/deceleration) to cross the intersection without any lateral or rear-end collision with the other CAVs, and without violating any of the state and control constraints.}
\subsection{Modeling Framework}
We model each CAV $i\in\mathcal{N}(t)$ as a double integrator
{
\begin{gather}%
\dot{p}_{i}(t)  =v_{i}(t),~
\dot{v}_{i}(t)  =u_{i}(t),~ { t\in [t_i^0, t_i^f],}
\label{eq:model2}%
\end{gather}}
where $p_{i}(t)\in\mathcal{P}_{i}$, $v_{i}(t)\in\mathcal{V}_{i}$, and $u_{i}(t)\in\mathcal{U}_{i}$ denote the position, speed and
acceleration (control input) of each CAV $i\in \mathcal{N}(t)$. 
The sets $\mathcal{P}_{i}$, $\mathcal{V}_{i}$, and $\mathcal{U}_{i}$, $i\in\mathcal{N}(t),$ are complete and totally bounded subsets of $\mathbb{R}$.
Let $\textbf{x}_{i}(t)=\left[p_{i}(t) ~ v_{i}(t)\right] ^{T}$ denote the state vector of each CAV $i$, with initial value $\textbf{x}_{i}^{0}=\left[p_{i}^{0} ~ v_{i}^{0}\right]  ^{T}$ taking values in $\mathcal{X}_{i}%
=\mathcal{P}_{i}\times\mathcal{V}_{i}$. The state space $\mathcal{X}_{i}$ for each CAV $i$ is closed with respect to the induced topology on $\mathcal{P}_{i}\times \mathcal{V}_{i}$ and thus, it is compact.

To ensure that the control input and speed of each CAV $i\in\mathcal{N}$ are within a
given admissible range, we impose the following constraints
\begin{gather}%
u_{i,min} \leq u_{i}(t)\leq u_{i,max},\quad\text{and} \nonumber \\ 
0  \leq v_{min}\leq v_{i}(t)\leq v_{max},\quad{ t\in\lbrack t_{i}%
^{0},t_{i}^{f}]},
\label{eq:state_control_constraint}%
\end{gather}
where $u_{i,min}$, $u_{i,max}$ are the minimum and maximum
acceleration for each CAV $i\in\mathcal{N}(t)$, and $v_{min}$, $v_{max}$ are the minimum and maximum speed limits respectively. Without loss of generality, we assume homogeneity in terms of CAV types, which enables the use of the same maximum acceleration $u_{max}$ and minimum acceleration $u_{min}$ for any CAV $i \in \mathcal{N}(t)$.
 %
{To ensure the avoidance} of rear-end collision of two consecutive CAVs traveling on the same lane, we impose the rear-end safety constraint 
%
\begin{equation}
s_{i}(t) \geq \delta_i(t),\quad \quad i \in \mathcal{N}(t),~ { t\in [t_i^0, t_i^f]},
\label{eq:rearend_constraint}
\end{equation}
where $s_{i}(t):=p_{k}(t)-p_{i}(t)$ is defined as the distance between CAV $i, k \in \mathcal{N}(t)$, where CAV $k$ is physically located immediately ahead of CAV $i$, and  $\delta_i(t)$ is the minimum safe distance which is a function of speed $v_i(t)$.
%
For each CAV $i\in \mathcal{N}(t)$, we define the set $\Gamma_i \colon= \{t\,\,|\,t\in  [t_i^{m}, t_i^f] \}$.
Lateral collision between any two CAVs $i,j\in \mathcal{N}(t)$ can be avoided if 
\begin{equation}
\Gamma_i \cap \Gamma_j=\varnothing, \quad \quad i,j\in \mathcal{N}(t),~ { t\in [t_i^{m}, t_i^f]}.
\label{eq:lateral_constraint}
\end{equation} 
In the modeling framework described above, we impose the following assumptions:
{
\begin{customass}{1} \label{assu:1}
Each CAV $i\in\mathcal{N}(t)$ communicates with each other and with the coordinator without any delays or errors.
\end{customass}}
 
\begin{customass}{2} \label{assu:2}
{For each CAV $i\in \mathcal{N}(t)$, no lane change maneuver is allowed  within the control zone.} 
\end{customass}

\begin{customass}{3} \label{assu:3}
None of the state constraints are active at time $t_i^0$ when each CAV $i \in \mathcal{N}(t)$ enters the control zone.
\end{customass}

The first assumption may be strong but it is relatively straightforward to relax it as long as the noise in the measurements and/or delays is bounded. For example, we can determine upper bounds on the state uncertainties as a result of sensing or communication errors and delays, and incorporate these into more conservative safety constraints. The second assumption allows us to focus only on the control of longitudinal vehicle dynamics {of CAVs within the control zone}. {Each CAV $i\in \mathcal{N}(t)$, however, can change lanes before the entry and/or after the exit of the control zone. Our analysis can include multiple lanes by appropriately revising the vehicle dynamics model \eqref{eq:model2}.} Finally, the third assumption ensures that, for each CAV $i\in\mathcal{N}(t)$, the initial state {at the entry of the control zone} is feasible. 

\subsection{Low-level Optimal Control Problem}
For each CAV $i\in\mathcal{N}(t)$, { $t\in [t_i^0, t_i^m]$}, traveling inside the control zone, we formulate the following optimal control problem
\begin{gather}\label{eq:decentral_problem}
\min_{u_i(t)\in U_i}  \int_{t_i^{0}}^{t_i^{m}} \frac{1}{2}u_i^2(t)~dt,\\ 
\text{subject to}%
:\eqref{eq:model2},\eqref{eq:state_control_constraint},\text{
}p_{i}(t_i^{0})=0\text{, }{p_{i}(t_i^{m})=L},\nonumber\\
\text{and given }t_i^{0}\text{, }v_{i}(t_i^{0})\text{, }t_i^{m},\nonumber
\end{gather}
where we consider the $L^2$-norm of the control input, i.e., {$u_i^2(t),$ as the cost function. By minimizing transient engine operation, we have direct benefits in fuel consumption in conventional vehicles (vehicles with internal combustion engines)}; see \cite{Malikopoulos2017}. 
Note that we do not explicitly include the lateral \eqref{eq:lateral_constraint} and rear-end \eqref{eq:rearend_constraint} safety constraints {in \eqref{eq:decentral_problem}}. The lateral collision constraint is enforced by selecting the appropriate merging time $t_i^m$ for each CAV $i$ in the upper-level throughput maximization problem. The activation of rear-end safety constraint can be avoided under certain conditions; see \cite{Malikopoulos2018c}. 

In our formulation, the state constraints are $\textbf{S}_i(t,\textbf{x}_i(t)) :=\left [ v_{i}(t) - v_{max} ~
v_{min} - v_{i}(t)\right ]^T\le 0$.
Note that, $\textbf{S}_i(t,\textbf{x}_i(t))$ is not an explicit function of the control input $u_i(t)$.
Thus, to formulate the tangency constraints, we need to take successive time derivatives of
$\textbf{S}_i(t,\textbf{x}_i(t))$ until we obtain an expression that is explicitly dependent on $u_i(t)$; {see \cite{bryson1975applied}}. If $q$ time derivatives are required, we refer to each constraint in $\textbf{S}_i^{(q)}(t,\textbf{x}_i(t))$ as the $q$th-order state variable inequality constraint. In our case, we have 1st-order speed constraint, e.g., $\textbf{S}_i^{(1)}(t,\textbf{x}_i(t),u_i(t)) = \left[
\begin{array}
[c]{ll}%
& \mbox{$ ~~u_i(t)$}\\
& \mbox{$ -u_i(t)$}
\end{array}
\right].$

To derive an analytical solution of the optimal control problem in \eqref{eq:decentral_problem} for each CAV $i\in\mathcal{N}(t)$, we formulate the adjoined Hamiltonian function $H_{i}\big(t, \textbf{x}_i(t),u_i(t)\big)$, { $t\in [t_i^0, t_i^m]$}, as follows,
\begin{align}\label{eq:lagrangian}
&H_{i}\big(t, \textbf{x}_i(t),u_i(t)\big) = \frac{1}{2} u^{2}_{i}(t) + \lambda^{p}_{i}(t) \cdot v_{i}(t) + \lambda^{v}_{i}(t) \cdot u_{i}(t)\nonumber\\ \nonumber 
&+ \boldsymbol{\mu}_i^T(t) \cdot \textbf{C}_i(t,\textbf{x}_i(t),u_i(t))+ \boldsymbol{\eta}_i^T(t) \cdot \textbf{S}_i(t,\textbf{x}_i(t)) \\ 
&=\frac{1}{2} u^{2}_{i}(t)
+ \lambda^{p}_{i}(t) \cdot
v_{i}(t)  + \lambda^{v}_{i}(t) \cdot u_{i}(t)  \\ \nonumber
&+\mu^{a}_{i}(t) \cdot(u_{i}(t) - u_{max})
+ \mu^{b}_{i}(t) \cdot(u_{min} - u_{i}(t)) \\ \nonumber
&+ \eta^{c}_{i}(t) \cdot(v_{i}(t)-v_{max})
+ \eta^{d}_{i}(t) \cdot(v_{min}-v_{i}(t)), %
\end{align}
where, {$\textbf{C}_i(t,\textbf{x}_i(t),u_i(t)):=[u_{i}(t) - u_{max}~~ u_{min} - u_{i}(t)]^T$ is the vector of control constraints in \eqref{eq:state_control_constraint}}, $\lambda^{p}_{i}(t),~  \lambda^{v}_{i}(t)$ are the co-state components corresponding to the state vector $\textbf{x}_i(t)$, and $\boldsymbol{\mu}_i(t)$ is the path co-vector for control constraints consisting of the Lagrange multipliers with the following conditions,
\begin{gather}
\mu^{a}_{i}(t) = \left\{
\begin{array}
[c]{ll}%
>0, & \mbox{$u_{i}(t) - u_{max} =0$},\\
=0, & \mbox{$u_{i}(t) - u_{max} <0$}, 
\end{array}
\right.\label{eq:kkt1}\\
\mu^{b}_{i}(t) = \left\{
\begin{array}
[c]{ll}%
>0, & \mbox{$u_{min} - u_{i}(t) =0$},\\
=0, & \mbox{$u_{min} - u_{i}(t)<0$},
\end{array}
\right.\label{eq:17b}
\end{gather}
and $\boldsymbol{\eta}_i(t)$ is the path co-vector for state constraints consisting of the Lagrange multipliers,
\begin{gather}
\eta^{c}_{i}(t) = \left\{
\begin{array}
[c]{ll}%
>0, & \mbox{$v_{i}(t) - v_{max} =0$},\\
=0, & \mbox{$v_{i}(t) - v_{max}<0$},
\end{array}
\right.\label{eq:17c} \\
\eta^{d}_{i}(t) = \left\{
\begin{array}
[c]{ll}%
>0, & \mbox{$v_{min} - v_{i}(t)=0$},\\
=0, & \mbox{$v_{min} - v_{i}(t)<0$}.
\end{array}
\right.\label{eq:16d}
\end{gather}

The {corresponding} Euler-Lagrange equations {at time $t\in [t_i^0, t_i^m]$} are
\begin{equation}\label{eq:EL1}
\dot\lambda^{p}_{i}(t) = - \frac{\partial H_i}{\partial p_{i}} = 0, \\
\end{equation}
\begin{equation}\label{eq:EL2}
\dot\lambda^{v}_{i}(t) = - \frac{\partial H_i}{\partial v_{i}} =\left\{
\begin{array}
[c]{ll}%
-\lambda^{p}_{i}(t) , & \mbox{$v_{i}(t) - v_{max} <0$}\\
\quad\quad\quad\quad\quad \text{and } &  \mbox{$v_{min} - v_{i}(t)<0$},\\
-\lambda^{p}_{i}(t)-\eta^{c}_{i}(t) , & \mbox{$v_{i}(t) - v_{max} =0$},\\
-\lambda^{p}_{i}(t)+\eta^{d}_{i}(t) , & \mbox{$v_{min} - v_{i}(t)=0$},
\end{array}
\right.
\end{equation}
and
{
\begin{equation}
\label{eq:EL3}\frac{\partial H_i}{\partial u_{i}} = u_{i}(t) + \lambda
^{v}_{i}(t) + \mu^{a}_{i}(t) - \mu^{b}_{i}(t) = 0.
\end{equation}}
If the inequality state and control constraints \eqref{eq:state_control_constraint} are not active, { we have $\mu^{a}_{i}(t) = \mu^{b}_{i}(t)= \eta^{c}_{i}(t)=\eta^{d}_{i}(t)=0$}. 
Applying the necessary conditions, the optimal control $u_i^*(t)$ can be derived from $u_i^*(t) + \lambda^{v}_{i}(t)= 0,~ i \in\mathcal{N}(t).$
From (\ref{eq:EL1}) and \eqref{eq:EL2} we have $\lambda^{p}_{i}(t) = a_{i}$, and $\lambda^{v}_{i}(t) = -\big(a_{i}\cdot t + b_{i}\big)$, where $a_{i}$ and $b_{i}$ are constants of integration corresponding to each CAV $i\in\mathcal{N}(t)$. Therefore, the {unconstrained} optimal control input $u_i^*(t)$ is
\begin{align}
&u^{*}_{i}(t) = a_{i} \cdot t + b_{i}, ~ { t\in [t_i^0, t_i^m]}. \label{eq:20}
\end{align}
Substituting the last equation into \eqref{eq:model2} we find the optimal speed and position for each CAV $i\in\mathcal{N}(t)$, namely
\begin{align}
&v^{*}_{i}(t) = \frac{1}{2} a_{i} \cdot t^2 + b_{i} \cdot t + c_{i}, \label{eq:21}\\
&p^{*}_{i}(t) = \frac{1}{6}  a_{i} \cdot t^3 +\frac{1}{2} b_{i} \cdot t^2 + c_{i}\cdot t + d_{i}, ~ { t\in [t_i^0, t_i^m]}, \label{eq:22}%
\end{align}
where $c_{i}$ and $d_{i}$ are constants of integration corresponding to each CAV $i \in \mathcal{N}(t)$. The constants of integration $a_i$, $b_{i}$, $c_{i}$, and $d_{i}$ can be determined from \eqref{eq:20}-\eqref{eq:22} using the initial and boundary conditions imposed in \eqref{eq:decentral_problem}. Note that, we can either compute $a_i$, $b_{i}$, $c_{i}$, and $d_{i}$ only once at time $t=t_i^0$ and apply the solution throughout optimization horizon $[t_i^0, t_i^m]$, or update the constants of integration by recomputing \eqref{eq:20}-\eqref{eq:22} at some discrete time step in $[t_i^0, t_i^m]$ to account for any disturbance within the control zone. For the remainder of the paper, we reserve the notations $a_i$, $b_{i}$, $c_{i}$, and $d_{i}$ only for the unconstrained optimal solution given in \eqref{eq:20}-\eqref{eq:22}.

\begin{customrem}{1}\label{rem:1}
For the case where the constants of integration $a_i = 0$ and $b_i=0$, we have the trivial solution of the unconstrained problem \eqref{eq:20}-\eqref{eq:22} as $u_i^*(t) =0,~v_i^*(t) = c_i, ~ p_i^*(t) = c_i\cdot t+d_i$, { $t\in [t_i^0, t_i^m]$}. This implies that if the speed is constant and the speed constraint is not active at time $t=t_i^0$ (Assumption \ref{assu:3}), none of the state and control constraints becomes active {for $t\in[t_i^0, t_i^m]$}. If $a_i,b_i\neq 0$, we have $u_i^*(t_i^0)\neq 0$.
\end{customrem}
In what follows, we only consider the non-trivial case {(Remark \ref{rem:1})} of the constrained optimization problem \eqref{eq:decentral_problem} where $a_i,b_i\neq 0$.

\section{Analysis of the Constrained Optimal Control Problem}

To derive the constrained analytical solution of \eqref{eq:decentral_problem}, we follow the standard methodology used in optimal control problems with interior point state and/or control constraints; see \cite{Bryson:1963,bryson1975applied}. Namely, we first start with the unconstrained arc and derive the solution using \eqref{eq:20}-\eqref{eq:22}. If the solution violates any of the state or control constraints, then the unconstrained arc is pieced together with the arc corresponding to the activated constraint, and we re-solve the problem with the two arcs pieced together at the junction point between the constrained and unconstrained arcs of the constrained solution \eqref{eq:decentral_problem}. The two arcs yield a set of algebraic equations which are solved simultaneously using the boundary conditions of \eqref{eq:decentral_problem} and the interior conditions between the arcs. If the resulting solution, which includes the determination of the junction point from one arc to the next one, violates another constraint, then the last two arcs are pieced together with the arc corresponding to the new activated constraint, and we re-solve the problem with the three arcs pieced together. The three arcs will yield a new set of algebraic equations that need to be solved simultaneously using the boundary conditions of \eqref{eq:decentral_problem} and interior conditions between the arcs. The resulting solution includes the junction point from one arc to the next one. The process is repeated until the solution does not violate any other constraints. 

This process can be computationally intensive for the following reasons. First, the recursive solution process to resolve all possible combinations of constraint activation might lead to intensive computation that prohibits real-time implementation. Second, each of the aforementioned recursion needs to be solved numerically due to the presence of implicit functions. To address both issues, we introduce a condition-based framework for the optimal control problem in \eqref{eq:decentral_problem} which leads to a closed-form analytical solution without this recursive procedure.


\subsection{Condition of Constraint Exclusion}\label{sec:step1}
For the optimal control problem in \eqref{eq:decentral_problem}, we have two state and two control constraints leading to 
$15$ possible constraint combinations in total that can become active within the optimization horizon $[t_i^0, t_i^m]$. In this section, we show that it is only possible for a subset of the constraints to become active in $[t_i^0, t_i^m]$. Therefore, it is not necessary to consider all the cases in \eqref{eq:decentral_problem}. In what follows, we delve deeper into the nature of the unconstrained optimal solution given in \eqref{eq:20}-\eqref{eq:22} to derive useful information about the possible existence of constraint activation within the control zone.
\begin{customlem}{1} \label{lem:1}
For each CAV $i\in\mathcal{N}(t)$, let $a_i$ and $b_i$ be the constants of integration of the unconstrained solution of \eqref{eq:decentral_problem} corresponding to the optimal control input $u_i^*(t)$, {$t\in [t_i^0, t_i^m]$}. If the speed $v_i(t)$ is not specified at $t_i^m$, then 
\begin{equation}
    a_i\cdot t_i^m+b_i=0, \quad t_i^m > t_i^0\ge0. \label{eq:b_i}
\end{equation}
\end{customlem} 
\begin{pf}
For all $i\in \mathcal{N}(t)$, since the speed $v_i(t)$ at $t=t_i^m > t_i^0$ is not fixed, we have $ \lambda_i^v(t_i^m) =  0$ \citep{naidu2002optimal}, which implies $u_i^*(t_i^m)=0$, and the result follows.
\end{pf} 
\begin{customcor}{1}\label{cor:1}
The constants of integration $a_i$ and $b_i$ of the unconstrained solution of \eqref{eq:decentral_problem} have opposite signs. 
\end{customcor}
\begin{pf}
Since $t_i^m$ is positive and non-zero, the result follows from \eqref{eq:b_i}.
\end{pf}
\begin{customcor}{2}\label{cor:2}
The {unconstrained} optimal control input $u_i^*(t)$ is  {linearly} either increasing or decreasing with respect to time, and $u_i^*(t_i^m)=0.$ 
\end{customcor}
\begin{pf}
From \eqref{eq:20}, $u_i^*(t)$ is a linear function with $u_i^*(t_i^0) \neq 0$ for the non-trivial case (Remark \ref{rem:1}), and $u_i^*(t_i^m)=0$ (Lemma \ref{lem:1}), so the result follows. 
\end{pf}
\begin{customrem}{2}\label{rem:2}
The constants of integration $a_i$ and $b_i$ of the unconstrained solution of \eqref{eq:decentral_problem} represents the slope of $u_i^*(t),~t\in[t_i^0,t_i^m]$, and the initial value of the control input $u_i^*(t)$ at time $t=t_i^0$, respectively. 
\end{customrem}
\begin{customlem}{2}\label{lem:2}
Let $v_i(t_i^0)$ be the initial speed of CAV $i\in\mathcal{N}(t)$ when it enters the control zone at  $p_i(t_i^0)$ and travels up to the entry of the merging zone at $p_i(t_i^m)$. Then the nature of the unconstrained optimal control input $u_i^*(t)$ can be characterized using the following conditions based on the boundary conditions of $v_i(t_i^0), p_i(t_i^0)$ and $p_i(t_i^m)$:
(i) The unconstrained optimal control input $u_i^*(t)$ is linearly decreasing if $v_i(t_i^0)<\frac{(p_i(t_i^m)-p_i(t_i^0))}{t_i^m}$.
(ii) The unconstrained optimal control input $u_i^*(t)$ is linearly increasing if $v_i(t_i^0)>\frac{(p_i(t_i^m)-p_i(t_i^0))}{t_i^m}$.
\end{customlem} 
\begin{pf}
From \eqref{eq:21} and \eqref{eq:22}, we can write $v_i(t_i^0)=\frac{1}{2}a_i\cdot(t_i^0)^{2}+b_i\cdot t_i^0+c_i$ and $p_{i}(t_i^0) = \frac{1}{6}  a_{i} \cdot (t_i^0)^3 +\frac{1}{2} b_{i} \cdot (t_i^0)^2 + c_{i}\cdot t_i^0 + d_{i}$.
Without loss of generality, if we let $t_i^0=0$, we have
\begin{gather}
    c_i=v_i(t_i^0),~ d_i = p_i(t_i^0). \label{eq:c_i}
\end{gather}
Evaluating \eqref{eq:22} at $t=t_i^m$, we have $p_i(t_i^m)=\frac{1}{6}a_i\cdot(t_i^m)^3+\frac{1}{2}b_i\cdot(t_i^m)^2+c_i\cdot t_i^m+d_i.$
Substituting \eqref{eq:b_i} and \eqref{eq:c_i} in the above equation and solving for $a_i$, we have
\begin{equation}
  a_i=\frac{3(v_i(t_i^0)\cdot t_i^m-(p_i(t_i^m)-p_i(t_i^0)))}{(t_i^m)^3}. \label{eq:a_i}
\end{equation}
Since $t_i^m>0$, we have a non-positive constant of integration $a_i$, if $(v_i(t_i^0) \cdot t_i^m-(p_i(t_i^m)-p_i(t_i^0)))<0$. From Corollary \ref{cor:2} and Remark \ref{rem:2}, a non-positive $a_i$ indicates a negative slope for $u_i^*(t)$, which implies that $u_i^*(t)$ is a linearly decreasing acceleration, and the proof is complete.
The second part of Lemma \ref{lem:2} can be proved following similar steps, hence it is omitted.
\end{pf}

\begin{customrem}{3}\label{rem:3}
When the CAV $i\in\mathcal{N}(t)$ travels with its initial speed {$v_i(t_i^0)$} throughout the control zone, we have {$v_i(t_i^0) \cdot t_i^m=(p_i(t_i^m)-p_i(t_i^0))$}. From \eqref{eq:a_i}, this implies that $a_i=0$, referring to an optimal control input $u_i^*(t)$ with horizontal slope. Since $u_i^*(t_i^m)=0$ (Lemma \ref{lem:1}), {we have} $u_i^*(t)=0$, for all {$t\in {[t_i^0, t_i^m]}$}.
\end{customrem}

\begin{customlem}{3} \label{lem:3}
For the unconstrained optimal solution of \eqref{eq:decentral_problem}, if either $v_i(t)-v_{max}\le0$ or $u_i(t)-u_{max}\le0$ becomes active at any time $t\in [t_i^0, t_i^m]$, neither $v_{min}-v_{i}(t)\le0$ nor $u_{min}-u_{i}(t)\le0$ can become active in $[t_i^0, t_i^m]$. The reverse also holds.
\end{customlem}

\begin{pf}
{Let $u_i^*(t)=a_i \cdot t + b_i> 0>u_{min}$ at some time $t \in [t_i^0,t_i^m)$. Since $u_i^*(t_i^m)=0$ (Lemma \ref{lem:1}) and $u_i^*(t)$ is a linearly decreasing function (Corollary \ref{cor:2}), we have $u_i^*(t)>u_{min}$, for all $t \in [t_i^0,t_i^m]$, i.e., the constraint $u_{min}-u_{i}(t)\le0$ can not become active at any time in $t \in [t_i^0,t_i^m]$.} 
The corresponding quadratic optimal speed profile $v_i^*(t)$ in \eqref{eq:21} is a parabolic function of degree 2 with y-symmetric axis located at $t_i^m$ in the speed-time graph. Applying the necessary and sufficient condition of optimality in \eqref{eq:21}, we have
\begin{gather}
   {\dot{v}_i^*(t)} = a_i\cdot t + b_i=0, ~
    {\ddot{v}_i^*(t)} = a_i,~ {t\in[t_i^0, t_i^m].} \label{eq:sufficient}
\end{gather}
Solving the first equation of \eqref{eq:sufficient}, we have the {extremum} point at $t=-\frac{b_i}{a_i}$ which corresponds to the vertex of the parabola of \eqref{eq:21} at $t=t_i^m$. Whether this point corresponds to the maximum or minimum of the \eqref{eq:21} can be determined from the second part of \eqref{eq:sufficient}. {Since} $u_i^*(t)$ is decreasing, $a_i<0$ {(Remark \ref{rem:2})}. Thus, the second equation of \eqref{eq:sufficient} indicates a maximum value at the vertex $t_i^m$, indicating a concave quadratic profile of $v_i^*(t)$. Since the {extremum of the quadratic profile of $v_i^*(t)$ is located at $t_i^m$ and $v_{min}<v_i(t_i^0)<v_{max}$ (Assumption \ref{assu:3}), we have $v_i^*(t)>v_{min}$ for all $t \in [t_i^0,t_i^m]$.}
Therefore, the constraints $v_{min}-v_{i}(t)\le 0$ can not become active at any time $t\in [t_i^0, t_i^m]$, and the proof of the first part of Lemma \ref{lem:3} is complete.

{Conversely, let $u_i^*(t)=a_i \cdot t + b_i< 0<u_{max}$ at some $t\in [t_i^0, t_i^m]$. Since $u_i^*(t_i^m)=0$ (Lemma \ref{lem:1}) and $u_i^*(t)$ is linearly increasing in $t\in [t_i^0, t_i^m]$ (Remark \ref{rem:2}), $u_{min}-u_{i}(t)\le0$ can not become active at any $t\in [t_i^0, t_i^m]$}. In addition, $u_i^*(t)$ yields a {convex quadratic profile of $v_i^*(t)$} with vertex at $t=t_i^m$.{ Since the extremum point is located at $t_i^m$ and  $v_{min}<v_i(t_i^0)<v_{max}$ (Assumption \ref{assu:3}), we have $v_i^*(t)<v_{max}$ for any $t \in [t_i^0,t_i^m]$}, which implies that the state constraint $v_i(t) - v_{max} \le 0$ cannot become active {at any time $t \in [t_i^0,t_i^m]$}.
\end{pf}

\begin{customcor}{3}\label{cor:3}
The sign of $a_i$ corresponding to the unconstrained solution of \eqref{eq:decentral_problem} dictates the activation of either constraint set $\{v_i(t)-v_{max}\le0$, $u_i(t)-u_{max}\le0\}$ or \{$v_{min}-v_{i}(t)\le0$, $u_{min}-u_{i}(t)\le0$\}.
\end{customcor}
\begin{pf}
Since $a_i$ is the slope of the optimal control input $u_i^*(t)$ (Remark \ref{rem:2}), the sign of $a_i$ determines whether $u_i^*(t)$ is positive and decreasing or negative and increasing, which, in turn, determines the constraint activation criteria in Lemma \ref{lem:3}.
\end{pf}

\begin{customrem}{4}\label{rem:4}
{The sign of $a_i$ can provide direct insight on which of the state and control constraints becomes active, and thus it can reduce the cardinality of the set of possible constrain activation cases.}
\end{customrem}

Based on Lemmas \ref{lem:2} and \ref{lem:3}, we now present the following result which provides the condition under which the state and control constraints  become active. { Note that the result is based on the initial and final conditions of \eqref{eq:decentral_problem}} which enable the determination of the possible constraint activation set without solving the unconstrained optimization problem in \eqref{eq:decentral_problem}.

\begin{customthm}{1}\label{theo:1}
Let CAV $i\in\mathcal{N}(t)$ enter the control zone with initial speed $v_i(t_i^0)$ and travel with the unconstrained optimal control input {$u_i^*(t)$, $t\in[t_i^0, t_i^m]$.} {Then, (i) $v_{min}-v_{i}(t)\le 0$ and $u_{min}-u_{i}(t)\le 0$ do not become active in $t\in[t_i^0, t_i^m]$, if $v_i(t_i^0)<\frac{(p_i(t_i^m)-p_i(t_i^0))}{t_i^m}$, and (ii)  $v_{i}(t)-v_{max}\le 0$ and $u_i(t)-u_{max}\le 0$ do not become active in $t\in[t_i^0, t_i^m]$, if $v_i(t_i^0)>\frac{(p_i(t_i^m)-p_i(t_i^0))}{t_i^m}$}. 
\end{customthm}

\begin{pf}
If $v_i(t_i^0)<\frac{(p_i(t_i^m)-p_i(t_i^0))}{t_i^m}$, then from \eqref{eq:a_i} $a_i<0$, hence $u_i^*(t)$ is linearly decreasing (Lemma \ref{lem:2}). Therefore, from Lemma \ref{lem:3}, $v_{min}-v_{i}(t)\le 0$ and $u_{min}-u_i(t)\le 0$ can not become active in $t\in[t_i^0, t_i^m]$, which concludes the proof of the first part. 

For the second part of Theorem 1, suppose that $v_i(t_i^0)>\frac{(p_i(t_i^m)-p_i(t_i^0))}{t_i^m}$. Hence $a_i>0$ (Lemma \ref{lem:2}), and $u_i^*(t)$ is linearly increasing. Therefore, from Lemma \ref{lem:3}, $v_i(t)-v_{max}(t)\le 0$ and $u_i(t)-u_{max}\le 0$ can not become active in {$t \in [t_i^0, t_i^m]$}, and the proof is complete.
\end{pf}

\begin{customrem}{5}\label{rem:5}
Theorem \ref{theo:1} aims at reducing the possible set of constraint activation cases. For example, if the condition in part (i) of Theorem \ref{theo:1} holds, then from the 15 possible cases of constraint activation, we only need to consider 3 cases: (a) $v_i(t)-v_{max}\le 0$, (b) $u_i(t)-u_{max}\le 0$, and (c) both  $v_i(t)-v_{max}\le 0$ and $u_i(t)-u_{max}\le 0$. Similarly, if the condition in part (ii) of Theorem 1 holds, then from the 15 possible cases of constraint activation, we only need to consider 3 cases: (a) $v_{min}-v_{i}(t)\le 0$, (b) $u_{min}-u_{i}(t)\le 0$, and (c) both $v_{min}-v_{i}(t)\le 0$ and $u_{min}-u_{i}(t)\le 0$.
\end{customrem}

Although Theorem \ref{theo:1} aims at reducing the possible constraint activation cases, it does not lead to the identification of the exact constraint activation of the unconstrained solution of \eqref{eq:decentral_problem}. In what follows, we provide the conditions that can be used to extend the results of Theorem \ref{theo:1} and identify the activation of any constraint case in $[t_i^0, t_i^m]$.
\subsection{Conditions of Constraint Activation}\label{sec:step2}

We start our exposition with some results that contain essential properties of the state and control constraint activation.

\begin{customlem}{4}\label{lem:4}
If neither $u_i(t)-u_{max}\le 0$ nor $u_{min}-u_{i}(t)\le 0$ is active at $t= t_i^0$, then it is guaranteed that neither of them will become active for all $t\in [t_i^0, t_i^m]$.
\end{customlem} 

\begin{pf}
Suppose that the unconstrained optimal solution of \eqref{eq:decentral_problem} yields  $u_i^*(t)=a_i t +b_i$ with  $a_i<0$. From Corollary \ref{cor:2} and Remark \ref{rem:2}, $u_i^*(t)$ decreases with respect to $t$, and at $t_i^m$, $u_i^*(t_i^m)=0$. Therefore, if $u_i^*(t_i^0) < u_{max}$, then $u_i^*(t)< u_{max}$ for all $t\in [t_i^0, t_i^m]$.  The second part of Lemma \ref{lem:4} can be proved following similar steps, hence it is omitted.
\end{pf}

\begin{customlem}{5}\label{lem:5}
If either $v_i(t)-v_{max}\le 0$ or $v_{min}-v_{i}(t)\le 0$ becomes active at any time $t\in [t_i^0, t_i^m)$, then it will remain active until $t=t_i^m$.
\end{customlem} 

\begin{pf}
Suppose that the unconstrained optimal solution of \eqref{eq:decentral_problem} yields  $u_i^*(t)=a_i t +b_i$ with  $a_i<0$. From Corollary \ref{cor:2} and Remark \ref{rem:2}, $u_i^*(t)$ decreases with respect to $t$, and at $t_i^m$, $u_i^*(t_i^m)=0$, which implies that $v_i^*(t)$ is monotonically increasing, i.e., $v_i^*(t_i^m) \ge v_i^*(t)$ {in $t\in [t_i^0, t_i^m)$}. Therefore,  $v_i(t)-v_{max}\le 0$ will remain active until $t=t_i^m.$
The second part of Lemma \ref{lem:5} can be proved following similar steps, hence it is omitted.
\end{pf}

\begin{customrem}{6}\label{rem:6}
Lemma \ref{lem:4} implies that the entry of the control-constrained arc can be only at $t= t_i^0,$ while Lemma \ref{lem:5} implies that there is no exit point in $[t_i^0, t_i^m]$ of the state-constrained arc after it becomes active.
\end{customrem}

The following results provide the conditions for which state and control constraint activation cases can be identified for the optimal control problem \eqref{eq:decentral_problem} a priori. 

\begin{customthm}{2}\label{theo:2}
Let $u_i^*(t)=a_i t +b_i$, {$t\in [t_i^0, t_i^m]$}, be the optimal control input of CAV $i\in\mathcal{N}(t)$ for the unconstrained solution of \eqref{eq:decentral_problem}. Then,
(i) for $a_i<0$, $v_i(t)-v_{max}\le 0$ becomes active if $t_i^m\le \frac{3(p_i(t_i^m)-p_i(t_i^0))}{v_i(t_i^0)+2v_{max}}$, and
(ii) for $a_i>0$, $v_{min}-v_{i}(t)\le 0$ becomes active if $t_i^m\ge\frac{3(p_i(t_i^m)-p_i(t_i^0))}{v_i(t_i^0)+2v_{min}}$.
\end{customthm}
\begin{pf}
For $a_i<0$, suppose that there exists a time $t_i^s \in (t_i^0,t_i^m]$ at which $v_i(t)-v_{max}\le 0$ becomes active. Then, from \eqref{eq:21} and \eqref{eq:c_i}, we have $\frac{1}{2}a_i\cdot (t_i^s)^2+b_i \cdot t_i^s+v_i(t_i^0)=v_{max}$.
Solving the quadratic equation for $t_i^s$, we have $t_i^s = \frac{-2b_i\pm \sqrt{4b_i^2-8a_i\cdot(v_i(t_i^0)-v_{max})}}{2a_i}$, which yields $t_i^{s}=t_i^m \pm \sqrt{\frac{4b_i^2-8a_i\cdot(v_i(t_i^0)-v_{max})}{4a_i^2}}$. Since $t_i^{s}\le t_i^m$, a feasible solution of $t_i^s$ exists
if we have $\sqrt{4b_i^2-8a_i\cdot(v_i(t_i^0)-v_{max})}\ge 0$ resulting in $a_i\le \frac{2(v_i(t_i^0)-v_{max})}{(t_i^m)^2}$. Combining with \eqref{eq:a_i}, the proof of the first statement of Theorem \ref{theo:2} follows. 

For $a_i>0$, suppose that there exists a time $t_i^s \in (t_i^0,t_i^m]$ at which the state constraint $v_{min}-v_{i}(t)\le 0$ becomes active. Then, from \eqref{eq:21} and \eqref{eq:c_i}, we have $\frac{1}{2}a_i\cdot (t_i^s)^2+b_i\cdot t_i^s+v_i(t_i^0)=v_{min}$. Solving the above equation for $t_i^s$, we have $t_i^s = \frac{-2b_i\pm \sqrt{4b_i^2-8a_i\cdot(v_i(t_i^0)-v_{min})}}{2a_i}$, which yields $t_i^{s}=t_i^m \pm \sqrt{\frac{4b_i^2-8a_i\cdot(v_i(t_i^0)-v_{min})}{4a_i^2}}$. 
Since $t_i^{s}\le t_i^m$, we need to have $\sqrt{4b_i^2-8a_i\cdot (v_i(t_i^0)-v_{min})}\ge 0$, and combining with \eqref{eq:a_i}, the proof of the second statement of Theorem \ref{theo:2} follows. 
\end{pf} 
\begin{customthm}{3}\label{theo:3}
Let $u_i^*(t)=a_i \cdot t +b_i$, {$t\in [t_i^0, t_i^m]$}, be the optimal control input of CAV $i\in\mathcal{N}(t)$ for the unconstrained solution of \eqref{eq:decentral_problem}. Then,
(i) for $a_i<0$, $u_i(t)-u_{max}\le 0$ becomes active if $t_i^m\le\frac{-3v_i(t_i^0)+\sqrt{9(v_i(t_i^0))^2+12u_{max}\cdot(p_i(t_i^m)-p_i(t_i^0))}}{2u_{max}}$, and
(ii) for $a_i>0$, $u_{min}-u_{i}(t)\le 0$ becomes active if $t_i^m\ge\frac{-3v_i(t_i^0)+\sqrt{9(v_i(t_i^0))^2+12 u_{min} \cdot (p_i(t_i^m)-p_i(t_i^0))}}{2 u_{min}}$.
\end{customthm}
\begin{pf}
For $a_i<0$, without loss of generality, we let $t_i^0=0$.
Given $v_i(t_i^0)$, $p_i(t_i^0)$ and $p_i(t_i^m)$, we will show that $t_i^m$ determines whether  $u_i(t)-u_{max}\le 0$ becomes active or not.
Let $\hat{t}_i^m$ be the value for which $u_i(t)-u_{max}\le 0$ becomes active at $t_i^0$, and $\hat{a}_i$, $\hat{b}_i$ the corresponding constants of integration.
Then from \eqref{eq:b_i} and \eqref{eq:a_i}, we can write {$\hat{b}_i= -\frac{3(v_i(t_i^0)\cdot\hat{t}_i^m-L)}{(\hat{t}_i^m)^2}=u_{max}$}, where $L=p_i(t_i^m)-p_i(t_i^0)=p_i(\hat{t}_i^m)-p_i(t_i^0)$, which can be reduced to {$u_{max}\cdot(\hat{t}_i^m)^2+3v_i(t_i^0)\cdot \hat{t}_i^m-3L=0$.}  The solution of the last equation yields
{$ \hat{t}_i^m=\frac{-3v_i(t_i^0)\pm\sqrt{9(v_i(t_i^0))^2+12u_{max}\cdot L)}}{2u_{max}}$.}
Since $\hat{t}_i^m>0$, {$\hat{t}_i^m=\frac{-3v_i(t_i^0)+\sqrt{9(v_i(t_i^0))^2+12u_{max}\cdot L)}}{2u_{max}}$}. Hence, for any $t_i^m$ such that $t_i^m\le \hat{t}_i^m$, $u_i(t)-u_{max}\le 0$ becomes active, and the proof of the first statement of Theorem \ref{theo:3} is complete.

For $a_i>0$, without loss of generality, we let $t_i^0=0$.
Let $\hat{t}_i^m$ be a value that $u_{min}-u_{i}(t)\le 0$ becomes active at $t_i^0$, and $\hat{a}_i$, $\hat{b}_i$ the corresponding constants of integration. Then from \eqref{eq:b_i} and \eqref{eq:a_i}, we can write {$\hat{b}_i= -\frac{3(v_i(t_i^0)\cdot\hat{t}_i^m-L)}{(\hat{t}_i^m)^2}=u_{min}$}, where $L=p_i(t_i^m)-p_i(t_i^0)=p_i(\hat{t}_i^m)-p_i(t_i^0)$, which can be reduced to {$u_{min}\cdot(\hat{t}_i^m)^2+3v_i(t_i^0)\cdot \hat{t}_i^m-3L=0$.}  The solution of the last equation yields
{$ \hat{t}_i^m=\frac{-3v_i(t_i^0)\pm\sqrt{9(v_i(t_i^0))^2+12u_{min}\cdot L)}}{2u_{min}}$}, from which we have the only admissible result {$\hat{t}_i^m=\frac{-3v_i(t_i^0)+\sqrt{9(v_i(t_i^0))^2+12u_{min}\cdot L)}}{2u_{min}}$.} {Hence, for any $t_i^m$ such that $t_i^m\ge \hat{t}_i^m$, $u_{min}-u_i(t)\le 0$ becomes active, and the proof of the second statement of Theorem \ref{theo:3} is complete.}
%
\end{pf}
\subsection{Interdependence of Constraint Activation Cases}
We have discussed so far the conditions under which any of the state and control constraints become active. Using these conditions, we can derive the analytical solution of \eqref{eq:decentral_problem}. 
However, the resulting solution might activate additional constrained arcs. Therefore, we need to be able to identify beforehand under which conditions any additional constrained arcs may become active. Next, we provide a set of conditions based on the junction point {where} transition between the constrained and unconstrained arcs occur.
\begin{customthm}{4}\label{theo:4}
For CAV $i\in\mathcal{N}(t)$, let $\tau_s^*\in(t_i^0,t_i^m]$ be the junction point of the state constrained arc {where} either $v_i(t)-v_{max}\le0$ or $v_{min}-v_i(t)\le0$ becomes active.
Then,
(i) $v_i(t)-v_{max}\le0$ may cause $u_i(t)-u_{max}\le0$ to become active,  if $  \tau_s^* \le \frac{-3v_i(t_i^0)+\sqrt{9(v_i(t_i^0))^2+12u_{max}\cdot(p_i^*(\tau_s^*)-p_i(t_i^0))}}{2u_{max}}$, and
(ii) $v_{min}-v_i(t)\le0$ may cause $u_{min}-u_i(t)\le0$ to become active, if $ \tau_s^* \ge \frac{-3v_i(t_i^0)+\sqrt{9(v_i(t_i^0))^2+12 u_{min} \cdot (p_i^*(\tau_s^*)-p_i(t_i^0))}}{2 u_{min}}$.
\end{customthm}

\begin{pf}
Suppose that $v_i(t)-v_{max}\le0$ becomes active at $\tau_s^*$, where $t_i^0< \tau_s^* \le t_i^m$. Then from \eqref{eq:model2}, $u_i^*(t)=0$ {in} $t\in [\tau_s^*, t_i^m]$ and $p_i(\tau_s^*) = p_i(t_i^m)-v_{max}\cdot (t_i^m-\tau_s^*)$. 
We will determine whether any control constraint $u_i(t)-u_{max}\le0$ becomes active in $t\in [t_i^0,\tau_s^*]$.
From Lemma \ref{lem:4}, the control constraint becomes active at $t=t_i^0$. 
Let $\hat{t}_i^m$ be the value that $u_i(t)-u_{max}\le 0$ becomes active at $t_i^0$, and $\hat{a}_i$, $\hat{b}_i$ the corresponding constants of integration. Without loss of generality, if we let $t_i^0=0$, then from \eqref{eq:b_i} and \eqref{eq:a_i} we can write, {$\hat{b}_i= -\frac{3(v_i(t_i^0)\cdot\hat{t}_i^m-(p_i^*(\tau_s^*)-p_i(t_i^0))}{(\hat{t}_i^m)^2}=u_{max}$},  where $p_i(\tau_s^*)-p_i(t_i^0)=p_i(\hat{t}_i^m)-p_i(t_i^0)$, which can be reduced to {$u_{max}\cdot (\hat{t}_i^m)^2+3v_i(t_i^0)\cdot \hat{t}_i^m-3(p_i^*(\tau_s^*)-p_i(t_i^0))=0$.}  The solution of the last equation yields
{$ \hat{t}_i^m=\frac{-3v_i(t_i^0) \pm \sqrt{9(v_i(t_i^0))^2+12u_{max}\cdot (p_i^*(\tau_s^*)-p_i(t_i^0)))}}{2u_{max}}$.}
Since $\hat{t}_i^m>0$, {$\hat{t}_i^m=\frac{-3v_i(t_i^0)+\sqrt{9(v_i(t_i^0))^2+12u_{max}\cdot (p_i^*(\tau_s^*)-p_i(t_i^0)))}}{2u_{max}}$.} Hence, for any  $\tau_s^*$ such that $\tau_s^*\le \hat{t}_i^m$, $u_i(t)-u_{max}\le 0$ becomes active, and the proof of the first statement of Theorem \ref{theo:4} is complete.

Suppose that $v_{min}-v_i(t)\le0$ becomes active at $\tau_s^*$, where $t_i^0< \tau_s^* \le t_i^m$. Then from \eqref{eq:model2}, $u_i^*(t)=0$ {in} $t\in [\tau_s^*, t_i^m]$ and $p_i(\tau_s^*) = p_i(t_i^m)-v_{min}\cdot (t_i^m-\tau_s^*)$. 
Let $\hat{t}_i^m$ be the value that $u_{min}-u_i(t)\le 0$ becomes active at $t_i^0$, and $\hat{a}_i$, $\hat{b}_i$ the corresponding constants of integration. Without loss of generality, if we let $t_i^0=0$, then from \eqref{eq:b_i} and \eqref{eq:a_i} we can write, {$\hat{b}_i= -\frac{3(v_i(t_i^0)\cdot\hat{t}_i^m-(p_i^*(\tau_s^*)-p_i(t_i^0))}{(\hat{t}_i^m)^2}=u_{min}$,}  where $p_i(\tau_s^*)-p_i(t_i^0)=p_i(\hat{t}_i^m)-p_i(t_i^0)$, which can be reduced to {$u_{min}\cdot (\hat{t}_i^m)^2+3v_i(t_i^0)\cdot \hat{t}_i^m-3(p_i^*(\tau_s^*)-p_i(t_i^0))=0$.}  The solution of the last equation yields {$\hat{t}_i^m=\frac{-3v_i(t_i^0) \pm \sqrt{9(v_i(t_i^0))^2 + 12u_{min} \cdot (p_i^*(\tau_s^*) - p_i(t_i^0)))}}{2u_{min}}$,} where
{$\hat{t}_i^m=\frac{-3v_i(t_i^0)+\sqrt{9(v_i(t_i^0))^2+12u_{min}\cdot (p_i^*(\tau_s^*)-p_i(t_i^0)))}}{2u_{min}}$ is the only admissible result.}
Hence, for any  $\tau_s^*$ such that $\tau_s^* \ge \hat{t}_i^m$, $u_{min}-u_i(t)\le 0$ becomes active, and the proof of the second statement of Theorem \ref{theo:4} is complete.
\end{pf}
\begin{customthm}{5}\label{theo:5}
For CAV $i\in\mathcal{N}(t)$, let $\tau_c^*\in(t_i^0,t_i^m]$ be the junction point of the control constrained arc where either $u_i(t)-u_{max}\le0$ or $u_{min}-u_i(t)\le0$ becomes active. Then,
(i) $u_i(t)-u_{max}\le0$ may cause $v_i(t)-v_{max}\le0$ to become active, if $ t_i^m\ge \tau_c^* - \frac{2(v_i(\tau_c^*)-v_{max})}{u_{max}}$, and
{(ii) $u_{min}-u_i(t)\le0$ may cause $v_{min}-v_i(t)\le0$ to become active, if $ t_i^m\ge \tau_c^* - \frac{2(v_i(\tau_c^*)-v_{min})}{u_{min}}$.}
\end{customthm}
\begin{pf}
Suppose that $u_i(t)-u_{max}\le0$ becomes active at $t_i^0$ (Remark \ref{rem:6}) with an exit time at $\tau_c^*\in(t_i^0,t_i^m]$.  Then from \eqref{eq:model2}, $u_i^*(t)=u_{max}$ in $t\in [t_i^0,\tau_c^*]$. Consequently, we have $v_i(\tau_c^*) = v_i(t_i^0)+u_{max}\cdot \tau_c^*$. 
We will determine whether any state constraint $v_i(t)-v_{max}\le0$ becomes active for the unconstrained arc within $t\in [\tau_c^*, t_i^m]$.
Suppose that there exists a time $t_i^s \in (\tau_c^*,t_i^m]$ at which $v_i(t)-v_{max}\le 0$ becomes active in $[\tau_c^*, t_i^m]$. Without loss of generality, if we let $\tau_c^* = 0$, then the constants of integration $\hat{a}_i, \hat{b}_i$ are given by $\hat{a}_i =- \frac{u_{max}}{\hat{t}_i^m}$ and $\hat{b}_i = u_{max}$ (Remark \ref{rem:2}), where $\hat{t}_i^m:= t_i^m - \tau_c^*$. From \eqref{eq:21} and \eqref{eq:c_i}, we have $\frac{1}{2}\hat{a}_i\cdot (t_i^s)^2+\hat{b}_i \cdot t_i^s+v_i(\tau_c^*)=v_{max}$.
Solving the quadratic equation for $t_i^s$, we have $t_i^s = \frac{-2\hat{b}_i\pm \sqrt{4\hat{b}_i^2-8\hat{a}_i\cdot(v_i(\tau_c^*)-v_{max})}}{2\hat{a}_i}$, which yields $t_i^{s}=\hat{t}_i^m \pm \sqrt{\frac{4\hat{b}_i^2-8\hat{a}_i\cdot(v_i(\tau_c^*)-v_{max})}{4\hat{a}_i^2}}.$ 
Since we require $t_i^{s}\le \hat{t}_i^m$, we need to have $\sqrt{4\hat{b}_i^2-8\hat{a}_i\cdot(v_i(\tau_c^*)-v_{max})}\ge 0$ resulting in $\hat{a}_i\le \frac{2(v_i(\tau_c^*)-v_{max})}{(\hat{t}_i^m)^2}$. By using the value of $\hat{a}_i$ in the above equation and simplifying, the proof of the first statement of Theorem \ref{theo:2} follows. 

For the second statement of Theorem \ref{theo:5}, suppose that there exists a time $t_i^s \in (\tau_c^*,t_i^m]$ at which $v_{min}-v_i(t)\le 0$ becomes active in $[\tau_c^*, t_i^m]$. Without loss of generality, if we let $\tau_c^* = 0$, then the constants of integration $\hat{a}_i, \hat{b}_i$ are given by $\hat{a}_i = -\frac{u_{min}}{\hat{t}_i^m}$ and $\hat{b}_i = u_{min}$ (Remark \ref{rem:2}), where $\hat{t}_i^m:= t_i^m - \tau_c^*$. From \eqref{eq:21} and \eqref{eq:c_i}, we have $\frac{1}{2}\hat{a}_i\cdot (t_i^s)^2+\hat{b}_i \cdot t_i^s+v_i(\tau_c^*)=v_{min}$.
Solving the quadratic equation for $t_i^s$, we have $t_i^s = \frac{-2\hat{b}_i\pm \sqrt{4\hat{b}_i^2-8\hat{a}_i\cdot(v_i(\tau_c^*)-v_{min})}}{2\hat{a}_i}$, which yields $t_i^{s}=\hat{t}_i^m \pm \sqrt{\frac{4\hat{b}_i^2-8\hat{a}_i\cdot(v_i(\tau_c^*)-v_{min})}{4\hat{a}_i^2}}.$ 
Since $t_i^{s}\le \hat{t}_i^m$, we need to have $\sqrt{4\hat{b}_i^2-8\hat{a}_i\cdot(v_i(\tau_c^*)-v_{min})}\ge 0$ resulting in $\hat{a}_i\le \frac{2(v_i(\tau_c^*)-v_{min})}{(\hat{t}_i^m)^2}$. By using the value of $\hat{a}_i$ in the above equation and simplifying, the proof of the second statement of Theorem \ref{theo:2} follows. 
\end{pf}
\begin{customrem}{7}\label{rem:7}
The conditions in Theorems 4 and 5 depend on the junction points $\tau_s^*$ and $\tau_c^*$ of the corresponding constraint activation cases, which can be derived analytically from the known boundary conditions of \eqref{eq:decentral_problem}. Since the derivation of such analytical solution requires additional information, we provide the analysis in the following section.
\end{customrem}
\section{Analytical Solution of the Constrained Optimal Control Problem} \label{sec:analytic-solution}
To derive the analytical solution of \eqref{eq:decentral_problem}, we present a condition-based framework consisting of the following steps. We first evaluate the condition stated in Theorem \ref{theo:1} to reduce the set of possible constraint activation cases (Remark \ref{rem:5}). Then using above result, we evaluate the conditions presented in Theorems \ref{theo:2} and \ref{theo:3} to determine whether any constraint has become active. If none of the constraints in \eqref{eq:state_control_constraint} becomes active, we simply derive the unconstrained solution using \eqref{eq:20}-\eqref{eq:22} and terminate the process. However, if the conditions in Theorems \ref{theo:2} and \ref{theo:3} indicate the activation of any constraint cases, we need to evaluate further the conditions in Theorems \ref{theo:4} and \ref{theo:5} to determine whether any additional constraints may become active within the constrained solution as a result of the constraint cases identified from Theorems \ref{theo:2} and \ref{theo:3}.
Once the nature of the final constraint activation case is identified using Theorems \ref{theo:4} and \ref{theo:5}, we then piece together the relevant unconstrained and constrained arcs that yield a set of algebraic equations which are solved simultaneously using the boundary conditions of {\eqref{eq:decentral_problem}} and interior conditions between the arcs.

Since we piece together multiple constrained and unconstrained arcs, we denote the constants of integration corresponding to each arc by $a_i^{(p)},b_i^{(p)},c_i^{(p)},d_i^{(p)}$, $p=1,2,\dots, N_{arc},$ where $N_{arc}\in\mathbb{N}$ is the total number of arcs pieced together in the constrained solution {and $p$} represents the position of the arcs in terms of their appearance in the optimal solution starting from $t_i^0$ to $t_i^m$. For $N_{arc}$ arcs, we have $(N_{arc}-1)$ junction points. At any junction point $\tau$, the states are continuous, namely,
\begin{gather}
    p_i(\tau^{-}) = p_i(\tau^{+}), ~
    v_i(\tau^{-}) = v_i(\tau^{+}), \label{eq:speed-continuity}
\end{gather} 
where, $\tau^{-}$ and $\tau^{+}$ represent the time instance right before and right after $\tau$, respectively.

In what follows, we present the closed form analytical solution of different cases of state and control constraint activation to derive the optimal input $u_i^*(t)$, {$t\in[t_i^0,t_i^m]$}, for each CAV $i\in\mathcal{N}(t)$. 
\begin{customcase}{1}
Only the state constraint $v_i(t)-v_{max}\le 0$ becomes active.
\end{customcase}
In this case, we have $\mu^{a}%
_{i}(t) = \mu^{b}_{i}(t)=\eta^{d}_{i}(t)=0$. From \eqref{eq:EL1}, \eqref{eq:EL2}, and \eqref{eq:EL3}, we have $u_i(t)+\lambda_i^v(t) =0, \label{eq:nece:case1}~
     \dot{\lambda}_i^p(t)=0, \text{ and }
     \dot{\lambda}_i^v(t)=-\lambda_i^p(t)-\eta_i^c(t). \label{eq:EL2:case1}$
By Lemma \ref{lem:5}, CAV $i\in\mathcal{N}(t)$ exits the constrained arc at $t=t_i^m$ which leads to a single junction point.
Let $\tau_s$, $t_i^0<\tau_s<t_i^m$, be the junction point and let $\tau_s^-$ and $\tau_s^+$ be the time instance just before and after time $\tau_s$. The optimal speed and control input on the constrained arc are
\begin{gather}
 v_i^*(t) = v_{max},~
 u_i^*(t) = 0,~ { t\in[\tau_s,t_i^m].} \label{eq:case1_u_opt}
\end{gather}
The jump conditions of the costates and the Hamiltonian at $\tau_s$ are
\begin{subequations}
\begin{align}
&{\lambda_i^p}(\tau_s^{-}) = {\lambda_i^p}(\tau_s^{+}) + \mathbf{\pi}_i\cdot\frac{\partial }{\partial p_i(t)}\begin{bmatrix}
v_i(t) - v_{max}
\end{bmatrix}\bigg{|}_{t=\tau_s}, 
\label{eq:jump_case1_1}\\
&{\lambda_i^v}(\tau^{-}_{s}) = {\lambda_i^v}(\tau^{+}_{s}) + \mathbf{\pi}_i\cdot\frac{\partial }{\partial v_i(t)}\begin{bmatrix}
v_i(t) - v_{max}
\end{bmatrix}\bigg{|}_{t=\tau_s},
\label{eq:jump_case1_2}\\
&H_i(\tau^{-}_{s}) = H_i(\tau^{+}_{s}) -  \mathbf{\pi}_i\cdot\frac{\partial }{\partial t}\begin{bmatrix}
v_i(t) - v_{max}
\end{bmatrix}\bigg{|}_{t=\tau_s}, \label{eq:jump_case1_3}
\end{align}
\end{subequations}
where $\pi_i$ is a constant Langrange multiplier determined so that $v_i(t)-v_{max}= 0$ is satisfied.
Note that, \eqref{eq:jump_case1_1}-\eqref{eq:jump_case1_3} imply possible discontinuity of the costates and the Hamiltonian at $t = \tau_s$. The state  variables are continuous at $t = \tau_s$. From \eqref{eq:jump_case1_3}, we have
\begin{gather}
     \frac{1}{2}u_i^2(\tau_s^-) + \lambda_i^p(\tau_s^-)\cdot v_i(\tau_s^-) 
    + \lambda_i^v(\tau_s^-)\cdot u_i(\tau_s^-)\nonumber\\
    +  \eta_i^c(\tau_s^-) \cdot (v_i(\tau_s^-)-v_{max})
    =  \frac{1}{2}u_i^2(\tau_s^+) + \lambda_i^p(\tau_s^+)\cdot v_i(\tau_s^+) \nonumber\\
    + \lambda_i^v(\tau_s^+)\cdot u_i(\tau_s^+)+ \eta_i^c(\tau_s^+)\cdot (v_i(\tau_s^+)-v_{max}).\label{eq:case1:hamiltonian}
\end{gather}
From the continuity of the states and since $v_i(\tau_s^+)=v_{max}$, $u_i(\tau_s^+)=0$, we have $\lambda_i^p(\tau_s^-)\cdot v_i(\tau_s^-) = \lambda_i^p(\tau_s^+)\cdot v_i(\tau_s^+)$. { The Lagrange multiplier $\eta_i^c(t)$ in \eqref{eq:17c}}, yields $\eta_i^c(\tau_s^-)\cdot (v_i(\tau_s^-)-v_{max})= \eta_i^c(\tau_s^+)\cdot (v_i(\tau_s^+)-v_{max}) =0 $. By combining the above equations, \eqref{eq:case1:hamiltonian} reduces to $\frac{1}{2}u_i^2(\tau_s^-)
    + \lambda_i^v(\tau_s^-)\cdot u_i(\tau_s^-) =0$, which implies that either $u_i(\tau_s^-) = 0$ or $\frac{1}{2}u_i(\tau_s^-) + \lambda_i^v(\tau_s^-) = 0$, or both. Since the second term contradicts $u_i(t)+\lambda_i^v(t) =0$, we have $u_i(\tau_s^-)=0$. 
The Lagrange multiplier $\eta_i^c(t)$ is $\eta_i^c(t) = \left\{
\begin{array}
[c]{ll}%
& \mbox{$0, \quad \quad \quad ~ \,\, if~v_i(t)<v_{max},~t\in[t_i^0,\tau_s)$},\\
& \mbox{$-\lambda_i^p(t), \quad if~v_i(t)=v_{max},~t\in[\tau_s,t_i^m]$}.\\
\end{array}
\right.$

Using the Euler-Lagrange equations, interior conditions, the initial and final boundary conditions, and the terminal condition of the costates, we can formulate a set of equations by piecing the unconstrained and constrained arcs together at time $t=\tau_s$. This results in a total number of $9$ equations that we need to solve simultaneously to compute $4+4+1=9$ variables corresponding to the constants of integration of unconstrained and constrained arc, and the junction point $\tau_s^*$ respectively. From \eqref{eq:20}-\eqref{eq:22} and the boundary conditions in \eqref{eq:decentral_problem}, we receive the following $4$ equations: $\frac{1}{2} a_i^{(1)}\cdot (t_i^0)^2 + b_i^{(1)} \cdot t_i^0 + c_i^{(1)} = v_i(t_i^0),
\frac{1}{6} a_i^{(1)} \cdot (t_i^0)^3 + \frac{1}{2} b_i^{(1)} \cdot (t_i^0)^2 + c_i^{(1)} \cdot t_i^0 + d_i^{(1)} = p_i(t_i^0), a_i^{(2)} \cdot t_i^m + b_i^{(2)} = 0, 
{\frac{1}{6} a_i^{(2)} \cdot (t_i^m)^3 + \frac{1}{2} b_i^{(2)} \cdot (t_i^m)^2 + c_i^{(2)} \cdot t_i^m + d_i^{(2)} = p_i(t_i^m)}.$ From the state and control continuity at the junction point $\tau_s$, we receive the remaining $5$ equations are,
{
\begin{subequations}
\begin{align}
&\frac{1}{2} a_i^{(1)}\cdot (\tau_s)^2 + b_i^{(1)} \cdot \tau_s + c_i^{(1)} = v_{max},\\
&a_i^{(1)} \cdot \tau_s + b_i^{(1)} = 0,\\
&\frac{1}{6} a_i^{(1)} \cdot (\tau_s)^3 + \frac{1}{2} b_i^{(1)} \cdot (\tau_s)^2 + c_i^{(1)} \cdot \tau_s + d_i^{(1)}\nonumber \\ 
&+ v_{max} \cdot (t_i^m - \tau_s) = p_i(t_i^m),\\
&\frac{1}{2} a_i^{(2)}\cdot (\tau_s)^2 + b_i^{(2)} \cdot \tau_s + c_i^{(2)} = v_{max},\\
&a_i^{(2)} \cdot \tau_s + b_i^{(2)} = 0,
\end{align}
\end{subequations}}
%
%
where $a_i^{(1)},b_i^{(1)},c_i^{(1)},d_i^{(1)}$ and $a_i^{(2)},b_i^{(2)},c_i^{(2)},d_i^{(2)}$ are the constants of integration for the unconstrained and constrained arcs, respectively. {The recursive process to solve the above set of equations cannot be computed in real time. Additionally, the computational speed and convergence of numerical methods are also sensitive to the initial guess of the variables, which impose additional burden on the real-time computation effort.} 
However, if the junction point $\tau_s^*$ can be derived as an explicit function of the initial and final boundary conditions, then the above set of equations can lead to a closed-form solution that can be solved analytically in real time. 

\begin{customlem}{6}\label{lem:6}
For CAV $i\in\mathcal{N}(t)$, let $\tau_s^*$ be the junction point between the unconstrained and constrained arc of the state constrained $v_i(t)-v_{max} \le 0$  solution. Then $\tau_s^*$ is an explicit function of $p_i(t_i^m), ~ v_{max}, ~ t_i^m$, and $v_i(t_i^0)$, and can be expressed as $ \tau_s^* = \frac{3 (p_i(t_i^m)-v_{max}\cdot t_i^m)}{(v_i(t_i^0)-v_{max})}$.
\end{customlem}
\begin{pf}
See Appendix \ref{app:3}.
\end{pf}

\begin{customcase}{2}
Only the control constraint $u_i(t)-u_{max}\le 0$ becomes active.
\end{customcase}
In this case, we have $\mu^{b}_{i}(t)= \eta^{c}_{i}(t)=\eta^{d}_{i}(t)=0$. From \eqref{eq:EL1}, \eqref{eq:EL2}, and \eqref{eq:EL3}, we have $u_i(t)+\lambda_i^v(t)+\mu_i^a(t) =0,\label{eq:nece:case2}~
     \dot{\lambda}_i^p(t)=0, \text{ and }
     \dot{\lambda}_i^v(t)=-\lambda_i^p(t).\label{eq:EL2:case2}$
By Lemma \ref{lem:4}, CAV $i\in \mathcal{N}(t)$ enters the constrained arc at time $t=t_i^0$ and has a single exit junction point. Let $\tau_c,~t_i^0<\tau_c<t_i^m$, be the junction point where the control constrained arc transitions into the unconstrained arc, and let $\tau_c^-$ and $\tau_c^+$ be the immediate left and the right instance of $\tau_c$. The optimal control input $u_i^*(t)$ at the junction point is $ u_i^*(\tau_c) = u_{max}$.
The jump conditions are $\lambda_i^p(\tau_c^-)-\lambda_i^p(\tau_c^+) = 0, \label{eq:jump_case2_1}
\lambda_i^v(\tau_c^-)-\lambda_i^v(\tau_c^+) = 0, \label{eq:jump_case2_2} \text{ and }
H_i(\tau_c^+)- H_i(\tau_c^-)=0, \label{eq:jump_case2_3}$
%
which imply continuity of the costates and the Hamiltonian at the junction point $t=\tau_c$. The last jump condition leads to $    \frac{1}{2}u_i^2(\tau_c^-) + \lambda_i^p(\tau_c^-)\cdot v_i(\tau_c^-)
    + \lambda_i^v(\tau_c^-)\cdot u_i(\tau_c^-)+ \nonumber
    \mu_i^a(\tau_c^-) \cdot (u_i(\tau_c^-)-u_{max}) =\frac{1}{2}u_i^2(\tau_c^+) + \lambda_i^p(\tau_s^+)\cdot v_i(\tau_c^+) \nonumber
    + \lambda_i^v(\tau_c^+)\cdot u_i(\tau_c^+)+ \mu_i^a(\tau_c^+)\cdot (u_i(\tau_c^+)-u_{max}) \label{eq:case2_continuity_u_1}$.
From the continuity of the state and costate $\lambda_i^p$ at $t=\tau_c$, we have $\lambda_i^p(\tau_c^-)\cdot v_i(\tau_c^-) = \lambda_i^p(\tau_c^+)\cdot v_i(\tau_c^+)$. Moreover,  \eqref{eq:kkt1} yields $ \mu_i^a(\tau_c^-) \cdot (u_i(\tau_c^-)-u_{max})=\mu_i^a(\tau_c^+)\cdot (u_i(\tau_c^+)-u_{max})=0$, which after simplification leads to either $u_i(\tau_c^+) = u_i(\tau_c^-)$ or $\frac{1}{2}(u_i(\tau_c^+) + u_i(\tau_c^-)) + \lambda_i^v(\tau_c^+)=0$, or both. Both equations lead to $u_i(\tau_c^+)=u_i(\tau_c^-)=u_{max}$.  
The Lagrange multiplier $\mu_i^a(t)$ is $\mu_i^a(t) = \left\{
\begin{array}
[c]{ll}%
& \mbox{$-\lambda_i^v(t) - u_{max}, ~   if ~t\in[t_i^0,\tau_c)$},\\
& \mbox{$0, \quad \quad \quad \quad \quad \quad  ~if~ t\in[\tau_c, t_i^m]$}.\\
\end{array}
\right.$

Using the Euler-Lagrange equations, jump conditions at the junction point, the initial and final boundary conditions, and the costate condition at $t=t_i^m$, we can formulate a set of equations by piecing the constrained and unconstrained arcs together at $t=\tau_c$. In this case, we have a constrained arc with constant parameters $a_i^{(1)},b_i^{(1)},c_i^{(1)},d_i^{(1)}$, followed by an unconstrained arc with constant parameters $a_i^{(2)},b_i^{(2)},c_i^{(2)},d_i^{(2)}$ pieced together at junction point $\tau_c$, leading to $4+4+1=9$ variables that need to be determined. At time $t = t_i^0$ and $t = \tau_c$, we have the following set of equations for the constrained arc,
\begin{subequations}
\begin{align}
    & a_i^{(1)}\cdot t_i^0 + b_i^{(1)} = u_{max}, \label{eq:case2_sol1} \\
    & a_i^{(1)}\cdot \tau_c + b_i^{(1)} = u_{max}, \label{eq:case2_sol2}\\
    & \frac{1}{2}a_i^{(1)}\cdot(t_i^0)^2 + b_i^{(1)}\cdot t_i^0 + c_i^{(1)} = v_i(t_i^0), \label{eq:case2_sol3}\\
    & \frac{1}{6} a_i^{(1)}\cdot (t_i^0)^3 + \frac{1}{2} b_i^{(1)}\cdot (t_i^0)^2 + c_i^{(1)} \cdot t_i^0 +  d_i^{(1)} = p_i(t_i^0). \label{eq:case2_sol4}
\end{align}
\end{subequations}
From \eqref{eq:case2_sol1} and \eqref{eq:case2_sol2}, considering $t_i^0=0$ without loss of generality, we have $a_i^{(1)}=0$ and $b_i^{(1)} = u_{max}$. Substituting in \eqref{eq:case2_sol3}, we have $c_i^{(1)} = v_i(t_i^0)$. Finally, solving \eqref{eq:case2_sol4}, $d_i^{(1)} = p_i(t_i^0)$. The following set of equations aim to determine the remaining constants of integration $a_i^{(2)},b_i^{(2)},c_i^{(2)},d_i^{(2)}$ of the exiting unconstrained arc and the junction point $\tau_c^*$
\begin{subequations}\label{eq:case2_analytic_sol}
\begin{align}
   & a_i^{(2)}\cdot\tau_c + b_i^{(2)} = u_{max}, \label{eq:case2_110}\\
  &  a_i^{(2)}\cdot t_i^m + b_i^{(2)} = 0, \label{eq:case2_111}\\
  &  \frac{1}{2}a_i^{(2)}\cdot\tau_c^2 + (b_i^{(2)}-u_{max})\cdot \tau_c + c_i^{(2)} - v_i^0 = 0, \label{eq:case2_112} \\
   & \frac{1}{6} a_i^{(2)}\cdot \tau_c^3 + \frac{1}{2} (b_i^{(2)}-u_{max})\cdot \tau_c^2 + ( c_i^{(2)} - v_i^0)\cdot \tau_c \nonumber \\
   & + {d_i^{(2)}-p_i(t_i^0) = 0}, \label{eq:case2_113}\\
  &  \frac{1}{6} a_i^{(2)}\cdot(t_i^m)^3 + \frac{1}{2}b_i^{(2)}\cdot(t_i^m)^2 + c_i^{(2)} \cdot t_i^m + d_i^{(2)} = p_i(t_i^m). \label{eq:case2_114}
\end{align}
\end{subequations}
\begin{customlem}{7}\label{lem:7}
For CAV $i\in\mathcal{N}(t)$, let $\tau_c^*$ be the junction point between the unconstrained and  control constraint $u_i(t)-u_{max}\le 0$  solution. Then $\tau_c^*$ can be expressed as an explicit function of $p_i(t_i^m),~ p_i(t_i^0), ~ u_{max}, ~ t_i^m$, and $v_i(t_i^0)$.
\end{customlem}
\begin{pf}
See Appendix \ref{app:4}.
\end{pf}
\begin{customcase}{3}
Both state constraint $v_i(t)-v_{max}\le 0$ and the control constraint $u_i(t)-u_{max}\le 0$ become active.
\end{customcase}
If both $u_i(t)-u_{max}\le 0$ and $v_i(t)-v_{max}\le 0$ become active, we derive the analytical solution combining the steps described in the previous two cases. In this case, we have $\mu^{b}_{i}(t)=\eta^{d}_{i}(t)=0$. From \eqref{eq:EL1}, \eqref{eq:EL2}, and \eqref{eq:EL3}, we have $u_i(t)+\lambda_i^v(t)+\mu_i^a(t) =0,~
     \dot{\lambda}_i^p(t)=0, \text{ and }
     \dot{\lambda}_i^v(t)=-\lambda_i^p(t)-\eta_i^c(t)$.
Let $\tau_c$ be the junction point that CAV $i\in\mathcal{N}(t)$ exits the control constrained arc and $\tau_s$ be the junction point that CAV $i$ enters the state constrained arc such that $t_i^0<\tau_c < \tau_s < t_i^m$. 
The optimal control input at the control constrained arc is $ u_i^*(t) = u_{max},$ for all $t\in[t_i^0, \tau_c]$. In the state constrained arc, we have $v_i^*(t)=v_{max},~u_i^*(t)=0,$ for all $t\in[\tau_s, t_i^m]$. From the jump conditions at the junction points $\tau_c$ and $\tau_s$, we have continuity in the state and control input. The Lagrange multipliers $\mu_i^a(t)$ and $\eta_i^c(t)$ are given by $\mu_i^a(t) = \left\{
\begin{array}
[c]{ll}%
& \mbox{$0, \quad \quad \quad \quad \quad \quad ~~~  t\in(\tau_c, t_i^m]$},\\
& \mbox{$-\lambda_i^v(t) - u_{max}, \quad  t\in[t_i^0, \tau_c]$},\\
\end{array}
\right\},$ 
and \\
$\eta_i^c(t) = \left\{
\begin{array}
[c]{ll}%
& \mbox{$0, \quad \quad \quad \quad \quad \quad  t\in[t_i^0, \tau_s)$},\\
& \mbox{$-\lambda_i^p(t), \quad \quad \quad ~~  t\in[\tau_s, t_i^m]$}.\\
\end{array}
\right\}.$

Solving \eqref{eq:case2_sol1}-\eqref{eq:case2_sol4}, considering $t_i^0=0$ without loss of generality, the constants of integration $a_i^{(1)},b_i^{(1)},c_i^{(1)},d_i^{(1)}$ of the control constrained arc are $a_i^{(1)}=0,~ b_i^{(1)} = u_{max}, ~ c_i^{(1)} = v_i^0$ and $d_i^{(1)} = p_i(t_i^0)$. The unconstrained arc with constants of integration $a_i^{(2)},b_i^{(2)},c_i^{(2)},\text{ and }d_i^{(2)}$ can consists of the following set of equations,
\begin{subequations}
     \begin{align}
    & a_i^{(2)}\cdot\tau_c + b_i^{(2)} = u_{max}, \label{eq:case3_119}\\
   & \frac{1}{2}a_i^{(2)}\cdot\tau_c^2 + (b_i^{(2)}-u_{max})\cdot \tau_c + c_i^{(2)} - v_i^0 = 0,\label{eq:case3_120} \\
   & \frac{1}{6} a_i^{(2)} \cdot\tau_c^3 + \frac{1}{2} (b_i^{(2)}-u_{max})\cdot \tau_c^2 + ( c_i^{(2)} - v_i^0) \cdot \tau_c \nonumber \\
   &+  (d_i^{(2)}-p_i(t_i^0)) = 0, \label{eq:case3_121}\\
   & a_i^{(2)} \cdot \tau_s + b_i^{(2)} = 0,\label{eq:case3_122} \\
    & \frac{1}{2}a_i^{(2)}\cdot\tau_s^2 + b_i^{(2)}\cdot\tau_s + c_i^{(2)} - v_{max} = 0, \label{eq:case3_126}\\
    & \frac{1}{6} a_i^{(2)} \cdot (\tau_s)^3 + \frac{1}{2} b_i^{(2)} \cdot (\tau_s)^2 + c_i^{(2)} \cdot \tau_s + d_i^{(2)} \nonumber\\
    & + v_{max} \cdot (t_i^m - \tau_s) = p_i(t_i^m).\label{eq:case3_128}
\end{align}
\end{subequations}
Finally, the state-constrained arc with constants of integration $a_i^{(3)},b_i^{(3)},c_i^{(3)},d_i^{(3)}$ consists of the following set of equations,
\begin{subequations}
\begin{align}
  & a_i^{(3)}\cdot t_i^m + b_i^{(3)} = 0,\label{eq:case3_124} \\
  &  a_i^{(3)}\cdot \tau_s + b_i^{(3)} = 0, \label{eq:case3_125}\\
   & \frac{1}{2}a_i^{(3)}\cdot\tau_s^2 - b_i^{(3)}\cdot \tau_s - c_i^{(3)} -v_{max}= 0,\label{eq:case3_127}\\
   & \frac{1}{6} a_i^{(3)}\cdot(t_i^m)^3 + \frac{1}{2}b_i^{(3)}\cdot(t_i^m)^2 + c_i^{(3)}\cdot t_i^m \nonumber  \\ 
  & + d_i^{(3)} - p_i(t_i^m) = 0\label{eq:case3_123}.
\end{align}
\end{subequations}
From \eqref{eq:case3_124}-\eqref{eq:case3_123}, we have $a_i^{(3)}=0,~ b_i^{(3)} = 0, ~ c_i^{(3)} = v_{max}$ and $d_i^{(3)} = p_i(t_i^m)-v_{max}\cdot t_i^m$.
The remaining constants of integration $a_i^{(2)},b_i^{(2)},c_i^{(2)},d_i^{(2)}$ of the unconstrained arc, and the junction points $\tau_s^*$ and $\tau_c^*$ can be determined by solving the set of equations \eqref{eq:case3_119}-\eqref{eq:case3_128}.

\begin{customlem}{8}\label{lem:8}
The junction point $\tau_s^*$ between the unconstrained and the constrained arc if $v_i(t)-v_{max} \le 0$ becomes active, and the junction point $\tau_c^*$ between the unconstrained and the constrained arc if $u_i(t)-u_{max} \le 0$ also becomes active are explicit functions of  $p_i(t_i^m), ~ v_{max},~u_{max}, ~ t_i^m$, and $v_i(t_i^0)$.
\end{customlem}
\begin{pf}
See Appendix \ref{app:5}.
\end{pf}

\begin{customcase}{4}
Only the state constraint $v_{min}-v_{i}(t)\le 0$ becomes active.
\end{customcase}
In this case, we have $\mu^{a}_{i}(t) = \mu^{b}_{i}(t)=\eta^{c}_{i}(t)=0$. From \eqref{eq:EL1}, \eqref{eq:EL2}, and \eqref{eq:EL3}, we have $u_i(t)+\lambda_i^v (t)=0, \label{eq:nece:case4}~
     \dot{\lambda}_i^p(t)=0, \text{ and }
     \dot{\lambda}_i^v(t)=-\lambda_i^p(t)-\eta_i^d(t). \label{eq:EL2:case4}$
Let $t = \tau_s$  be the junction point that $v_{min}-v_{i}(t)\le 0$ becomes active. The optimal speed and control at the junction point are $ v_i^*(t) = v_{min}, ~
 u_i^*(t) = 0,$ for all $t\in[\tau_s, t_i^m].$
%
%
The jump conditions are
\begin{subequations}
\begin{align}
&{\lambda_i^p}(\tau_s^{-}) = {\lambda_i^p}(\tau_s^{+}) + \mathbf{\pi}_i\cdot\frac{\partial }{\partial p_i(t)}\begin{bmatrix}
v_{min} - v_{i}(t)
\end{bmatrix}\bigg{|}_{t=\tau_s},
\label{eq:jump_case4_1}\\
&{\lambda_i^v}(\tau^{-}_{s}) = {\lambda_i^v}(\tau^{+}_{s}) + \mathbf{\pi}_i\cdot\frac{\partial }{\partial v_i(t)}\begin{bmatrix}
v_{min} - v_{i}(t)
\end{bmatrix}\bigg{|}_{t=\tau_s}, 
\label{eq:jump_case4_2}\\
&H_i(\tau^{-}_{s}) = H_i(\tau^{+}_{s}) -  \mathbf{\pi}_i\cdot\frac{\partial }{\partial t}\begin{bmatrix}
v_{min} - v_{i}(t)
\end{bmatrix}\bigg{|}_{t=\tau_s}, \label{eq:jump_case4_3}
\end{align}
\end{subequations}
where $\pi_i$ is a constant Langrange multiplier determined so that $v_{min}-v_{i}(t) = 0$ is satisfied.
Note that, \eqref{eq:jump_case4_1}-\eqref{eq:jump_case4_3} imply possible discontinuity of the costates and the Hamiltonian at $t = \tau_s$. The state  variables are continuous at $t = \tau_s$. From \eqref{eq:jump_case4_1} and \eqref{eq:jump_case4_3}, the position costate and the Lagrangian of the Hamiltonian is continuous at $t=\tau_s$. 

\begin{customlem}{9}\label{lem:9}
If the state constraint $v_{min}-v_{i}(t)\le 0$ becomes active, then the control input $u_i(t)$ is continuous at the junction point $t = \tau_s$.
\end{customlem}
\begin{pf}
See Appendix \ref{app:1}.
\end{pf}

The Lagrange multiplier $\eta_i^d(t)$ can be expressed as,
$\eta_i^d(t) =\bigg\{ \begin{matrix}
0, \quad \quad \quad \,\, if\quad t\in[t_i^0,\tau_s), \\
-\lambda_i^p(t), \quad if\quad ~t\in[\tau_s,t_i^m].
\end{matrix}$
Using the Euler-Lagrange equations, interior conditions, initial and final boundary conditions, and the costate condition at $t=t_i^m$, we can formulate a set of equations similar to Case 1 to solve for $4+4+1=9$ variables corresponding to the constants of integration of the unconstrained and constrained arc, and the junction point $\tau_s$. The set of equations of the unconstrained arc with constants of integration $a_i^{(1)},b_i^{(1)},c_i^{(1)},d_i^{(1)}$ are, $\frac{1}{2} a_i^{(1)}\cdot (t_i^0)^2 + b_i^{(1)} \cdot t_i^0 + c_i^{(1)} = v_i(t_i^0), ~
\frac{1}{6} a_i^{(1)} \cdot (t_i^0)^3 + \frac{1}{2} b_i^{(1)} \cdot (t_i^0)^2 + c_i^{(1)} \cdot t_i^0 + d_i^{(1)} = p_i(t_i^0), ~
\frac{1}{2} a_i^{(1)}\cdot (\tau_s)^2 + b_i^{(1)} \cdot \tau_s + c_i^{(1)} = v_{min}, ~
a_i^{(1)} \cdot \tau_s + b_i^{(1)} = 0, \text{ and }
\frac{1}{6} a_i^{(1)} \cdot (\tau_s)^3 + \frac{1}{2} b_i^{(1)} \cdot (\tau_s)^2 + c_i^{(1)} \cdot \tau_s + d_i^{(1)} + v_{min} \cdot (t_i^m - \tau_s) = p_i(t_i^m).$
The set of equations of the state constrained arc with the constants of integration $a_i^{(2)},b_i^{(2)},c_i^{(2)},d_i^{(2)}$ are $\frac{1}{2} a_i^{(2)}\cdot (\tau_s)^2 + b_i^{(2)} \cdot \tau_s + c_i^{(2)} = v_{min}, ~
a_i^{(2)} \cdot t_i^m + b_i^{(2)} = 0, ~
a_i^{(2)} \cdot \tau_s + b_i^{(2)} = 0, \text{ and }
\frac{1}{6} a_i^{(2)} \cdot (t_i^0)^3 + \frac{1}{2} b_i^{(2)} \cdot (t_i^0)^2 + c_i^{(2)} \cdot t_i^0 + d_i^{(2)} = p_i(t_i^m),$
which yield $a_i^{(2)}=0,~ b_i^{(2)} = 0, ~ c_i^{(2)} = v_{min}$ and $d_i^{(2)} = p_i(t_i^m)-v_{min}\cdot t_i^m$.
The remaining constants of integration $a_i^{(1)},b_i^{(1)},c_i^{(1)},d_i^{(1)}$ and the junction point $\tau_s^*$ can be determined numerically by solving simultaneously the above set of equations.

\begin{customlem}{10}\label{lem:10}
For CAV $i\in\mathcal{N}(t)$, let $\tau_s^*$ be the junction point between the unconstrained and constrained arc of the state constrained $v_{min}-v_{i}(t) \le 0$  solution. Then $\tau_s^*$ is an explicit function of $p_i(t_i^m), ~ v_{min}, ~ t_i^m$ and $v_i(t_i^0)$, and can be expressed as $ \tau_s^* = \frac{3 (p_i(t_i^m)-v_{min}\cdot t_i^m)}{(v_i(t_i^0)-v_{min})}$
\end{customlem}
\begin{pf}
The proof is similar to the proof of Lemma \ref{lem:6} (see Appendix \ref{app:3}), hence it is omitted.
\end{pf}

\begin{customcase}{5}
Only the control constraint $u_{min}-u_{i}(t)\le 0$ becomes active.
\end{customcase}
In this case, we have $\mu^{a}_{i}(t)= \eta^{c}_{i}(t)=\eta^{d}_{i}(t)=0$. From \eqref{eq:EL1}, \eqref{eq:EL2}, and \eqref{eq:EL3}, we have $u_i(t)+\lambda_i^v(t)-\mu_i^b(t) =0,\label{eq:nece:case5}~
     \dot{\lambda}_i^p(t)=0, \text{ and }
 \dot{\lambda}_i^v(t)=-\lambda_i^p(t).\label{eq:EL2:case5}$
Let $\tau_c>t_i^0$ be the junction point that CAV $i\in\mathcal{N}(t)$ transitions from the constrained arc to the unconstrained arc.  The optimal control at the junction point $\tau_c$ is $u_i^*(\tau_c) = u_{min}$.
From the jump conditions, we have $\lambda_i^p(\tau_c^-)=\lambda_i^p(\tau_c^+),$ 
$\lambda_i^v(\tau_c^-)=\lambda_i^v(\tau_c^+),$ and
$H_i(\tau_c^+)= H_i(\tau_c^-).$
\begin{customlem}{11}\label{lem:11}
If the control constraint $u_{min}-u_{i}(t)\le 0$ becomes active, then the control input $u(t)$ is continuous at the junction point $t = \tau_c$.
\end{customlem}
\begin{pf}
See Appendix \ref{app:2}.
\end{pf}
The Lagrange multiplier $\mu_i^b(t)$ can be expressed as,
$\mu_i^b(t) =\bigg\{ \begin{matrix}
\lambda_i^v(t) + u_{min}, \quad \quad if\quad t\in[t_i^0,\tau_c),\\
0, \quad \quad \quad \quad \quad \quad \quad if\quad t\in [\tau_c, t_i^m].

\end{matrix}$
Using the Euler-Lagrange equations, interior condition, initial and final boundary conditions, and the condition of costates at $t=t_i^m$, we have a set of equations of the constrained arc: $    a_i^{(1)}\cdot t_i^0 + b_i^{(1)} = u_{min}, \label{eq:case5_1} 
a_i^{(1)}\cdot \tau_c + b_i^{(1)} = u_{min}, \label{eq:case5_2}
    \frac{1}{2}a_i^{(1)}(t_i^0)^2 + b_i^{(1)}\cdot t_i^0 + c_i^{(1)} = v_i(t_i^0), \label{eq:case5_3} \text{ and }
    \frac{1}{6} a_i^{(1)} (t_i^0)^3 + \frac{1}{2} b_i^{(1)}\cdot (t_i^0)^2 + c_i^{(1)} \cdot t_i^0 +  d_i^{(1)} = 0. \label{eq:case5_4}$, resolving which with $t_i^0 = 0$ yields, $a_i^{(1)}=0,~b_i^{(1)}=u_{min},~c_i^{(1)}=v_i(t_i^0),~d_i^{(1)}=p_i(t_i^0)$,
where $a_i^{(1)},~b_i^{(1)},~c_i^{(1)},~d_i^{(1)}$ are the constants of integration for the constrained arc.
In addition, we have a set of equations of the unconstrained arc:
    $a_i^{(2)}\cdot\tau_c - b_i^{(2)} + u_{min} = 0, \label{eq:case5_5}
    a_i^{(2)}\cdot t_i^m + b_i^{(2)} = 0,  \label{eq:case5_6}
    \frac{1}{2}a_i^{(2)}\cdot \tau_c^2 + (b_i^{(2)}-u_{min})\cdot \tau_c + c_i^{(2)} - v_i(t_i^0) = 0,  \label{eq:case5_7}\\
    \frac{1}{6} a_i^{(2)}\cdot \tau_c^3 + \frac{1}{2}(b_i^{(2)}-u_{min})\cdot \tau_c + ( c_i^{(2)} - v_i(t_i^0)\cdot \tau_c +  d_i^{(2)} = 0,  \label{eq:case5_8} \text{ and }
    \frac{1}{6} a_i^{(2)} \cdot (t_i^m)^3 + \frac{1}{2}b_i^{(2)}\cdot (t_i^m)^2 + c_i^{(2)}\cdot t_i^m + d_i^{(2)} - p_i(t_i^m) = 0, \label{eq:case5_9}$
where $a_i^{(2)},~b_i^{(2)},~c_i^{(2)},~d_i^{(2)}$ are the constants of integration of the unconstrained arc.

\begin{customlem}{12}\label{lem:12}
For CAV $i\in\mathcal{N}(t)$, let $\tau_s^*$ be the junction point between the unconstrained and constrained arc of the control constrained ($u_{min}-u_{i}(t) \le 0$) solution of \eqref{eq:decentral_problem}. Then $\tau_c^*$ can be expressed as an explicit function of $p_i(t_i^m),~ p_i(t_i^0), ~ u_{min}, ~ t_i^m$, and $v_i(t_i^0)$.
\end{customlem}
\begin{pf}
The proof is similar to the proof of Lemma \ref{lem:7} (see Appendix \ref{app:4}), hence it is omitted.
\end{pf}

\begin{customcase}{6}
Both state constraint $v_{min}-v_{i}(t)\le 0$ and the control constraint $u_{min}-u_{i}(t)\le 0$ become active.
\end{customcase}
In this case, we can derive the analytical solution following similar steps to Case 3. 
A control constrained $u_{min}-u_{i}(t)\le 0$ arc with constants of integration $a_i^{(1)},b_i^{(1)},c_i^{(1)},d_i^{(1)}$ is pieced together with an unconstrained arc with constants of integration $a_i^{(2)},b_i^{(2)},c_i^{(2)},d_i^{(2)}$ at the junction point $\tau_c$. The unconstrained arc is pieced together with the state constrained $v_{min}-v_{i}(t)\le 0$ arc  with constants of integration $a_i^{(3)},b_i^{(3)},c_i^{(3)},d_i^{(3)}$ at the junction point $\tau_s$.
The constants of integration of the constrained and unconstrained arcs, and the junction points $\tau_s^*$ and $\tau_c^*$ can be determined by a set of equations similar to those derived in Case 3. 

\begin{customlem}{13}
The junction point $\tau_s^*$ between the unconstrained and the constrained arc when $v_{min}-v_{i}(t) \le 0$ becomes active, and the junction point $\tau_c^*$ between the unconstrained and the constrained arc when $u_{min}-u_{i}(t) \le 0$ also becomes active are explicit functions of $p_i(t_i^m), ~ v_{min},~u_{min}, ~ t_i^m$, and $v_i(t_i^0)$.
\end{customlem}
\begin{pf}
The proof is similar to the proof of Lemma \ref{lem:8} (see Appendix E), hence it is omitted.
\end{pf}

\section{Simulation Results}
We validate the analytical solution of the optimal control problem  \eqref{eq:decentral_problem} through numerical simulation in MATLAB. {In this section, we present the results considering $t_i^m = 10$ {s}, where only the state constraint $v_i(t)-v_{max}\le 0$ and control constraint $u_i(t)-u_{max}\le 0$ can become active (Theorem \ref{theo:1}). Similar results to those presented here can be also derived for the case where  $v_{min}-v_i(t)\le0$ and $u_{min}-u_i(t)\le0$ become active.} We consider the initial and final position {of CAV $i\in \mathcal{N}(t)$} to be $p_i(t_i^0)=0$ m and $p_i(t_i^m)=200$ m, and the {initial speed $v_i(t_i^0)=14.3\,$ m/s}. 
For each CAV $i\in\mathcal{N}(t)$, we enforce the maximum speed limit and acceleration to be $v_{max}=22$ {m/s} and $u_{max}=1.8$ m/s$^2$ respectively. 
\begin{figure}[ht]
\centering
\includegraphics[width=3.5in]{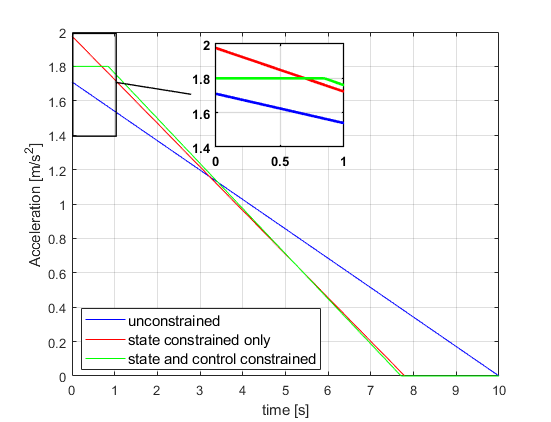} \caption{Optimal control trajectory for the unconstrained (blue), state constraint $v_i(t)-v_{max}\le 0$ only (red) and both state-control constraint (green) case.}%
\label{fig:2}%
\end{figure}
The standard procedure to solve the optimal control problem \eqref{eq:decentral_problem} is to identify whether any of the state or control constraints become active and derive the constrained solution in a recursive manner until none of the constraints are active, as shown in Fig. \ref{fig:2}. The unconstrained solution (blue trajectory in Fig. \ref{fig:2}) activates the state constraint $v_i(t)-v_{max}\le 0$ only. The acceleration corresponding to the state-constrained ($v_i(t)-v_{max}\le 0$) solution is shown by the red trajectory in Fig. \ref{fig:2}, where the unconstrained and constrained arcs are pieced together at the junction point at $t= 7.79$ {s}. However, the state-constrained solution (red trajectory in Fig. \ref{fig:2}) has to be re-derived since the control constraint $u_i(t)-u_{max}\le0$, which was not active before, becomes active now as shown by the red trajectory in Fig. \ref{fig:2}. The constrained optimal control input is derived by piecing the state and control constrained arcs together, and it is shown by the green trajectory in Fig. \ref{fig:2}.

In our condition-based framework, we do not need to consider the intermediate iterative steps above, i.e., the unconstrained (blue trajectory) and state constrained solution (red trajectory) in Fig. \ref{fig:2}. We can directly derive the final closed-form analytical solution (green trajectory in Fig. \ref{fig:2}) by sequentially checking the conditions in Theorems \ref{theo:1}-\ref{theo:5}. First, we start with Theorem \ref{theo:1} to reduce the possible constraint activation set. Since the first statement of Theorem \ref{theo:1} holds for $t_i^m=10$ {s} and the boundary conditions, we only need to consider whether $v_i(t)-v_{max}\le 0$ or $u_i(t)-u_{max}\le 0$ become active, which reduces the possible constraint activation cases from 15 to 3. 
Then, we use Theorems \ref{theo:2} and \ref{theo:3} to identify the specific constraint activation case. In this case, part (i) of Theorem \ref{theo:2} holds, indicating that $v_i(t)-v_{max}\le 0$ becomes active in $(t_i^0, t_i^m]$. However, part (i) of Theorem \ref{theo:3} does not hold indicating that $u_i(t)-u_{max}\le 0$ will not become active.
Using the result obtained above, we then check part (i) of Theorem \ref{theo:4} which readily indicates that an additional and initially non-existent control constraint {$u_i(t)-u_{max}\le 0$} becomes active within the state-constrained solution, as shown by the red trajectory in Fig. \ref{fig:2}. Using the result of Theorem \ref{theo:4}, we apply the analysis presented in Case 3 to determine the complete state and control constrained-optimal solution. {  Here, the aforementioned condition-based framework requires $0.001107$ s to solve in an Intel Core i7-6700 CPU @ 3.40 GHz using MATLAB R2017b.}
Note that, if the first statement of Theorem \ref{theo:4} does not hold, then none of the control constraints can become active, and thus we can use the analysis presented in Case 1 to determine the optimal solution.

\begin{figure}[ht]
\centering
\includegraphics[width=3.5in]{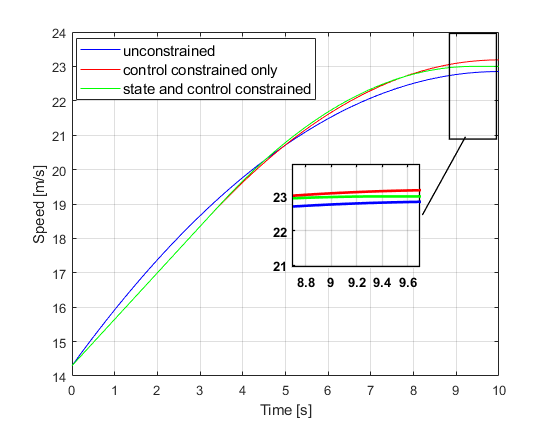} \caption{Optimal speed trajectory for the unconstrained (blue), control constraint $u_i(t)-u_{max}\le 0$ only (red) and both state-control constraint (green) case.}%
\label{fig:3}%
\end{figure}

Next, we consider a different scenario to show the impact when the control constraint ($u_i(t)-u_{max}\le0$) becomes active (Fig. \ref{fig:3}). In this case, we set the maximum speed $v_{max}$ and acceleration $u_{max}$ to be $23$ {m/s} and $1.35$ {m/s$^2$} respectively. Following the above procedure, we check part (i) of Theorem $\ref{theo:2}$ and \ref{theo:3}. Since only part (i) of Theorem \ref{theo:3} holds, we conclude that the control constraint $u_i(t)-u_{max}\le0$ will become active. We then check part (i) of Theorem \ref{theo:5} to check whether any additional state constraint will become active within the control constrained solution. In this case, part (i) of Theorem \ref{theo:5} holds, as evident from the control constrained state trajectory (red trajectory) in Fig. \ref{fig:3}. Therefore, we use the analysis presented in Case 3 to derive the complete state- and control- constrained solution as illustrated by the green trajectory in Fig. \ref{fig:3}. Note that, in Fig. \ref{fig:3}, in the unconstrained solution (blue trajectory) none of the state constraints become active. However, the control-constrained solution (red trajectory) activates the state constraint $v_i(t)-v_{max}\le0$. Based on our condition-based framework, we can avoid the computation of the intermediate solutions, i.e., the unconstrained trajectory (blue trajectory in Fig. \ref{fig:3}) and the control constrained trajectory (red trajectory in Fig. \ref{fig:3}), and directly derive the final constrained trajectory as illustrated by the green trajectory in Fig. \ref{fig:3}. 

\section{Concluding Remarks}
In this paper, we addressed the state and control constrained optimal framework for coordinating CAVs at different traffic scenarios {such as merging at roadways and roundabouts, cruising in congested traffic, passing through speed reduction zones, under 100\% CAV penetration}, and provided a condition-based framework to determine the constrained solution without requiring to follow the standard recursive process. We mathematically characterized the activation cases of different state and control constraint combinations, and provided  a priori conditions under which different constraint combination can become active. In addition, we presented the closed-form analytical solution of the constrained optimal control problem that can be derived and implemented in real time. We validated a subset of constraint activation cases through numerical simulation and showed how the proposed framework can identify the interdependent constraint activation based on the boundary conditions. By eliminating the intermediate steps of solving the constrained optimal control problem, the proposed condition-based framework improves on the standard methodology to solve the constrained optimal control problem.

The proposed framework has certain limitations since it does not consider the optimal control problem with constrained terminal speed, which may result in multiple junction points leading to a more complex formulation. {Moreover, in our framework, we considered 100\% penetration rate of CAVs having access to perfect information (no
errors or delays) which both impose limitations for real-world applications. It is expected that CAVs will gradually penetrate the market, interact with non-CAVs
and contend with vehicle-to-vehicle and vehicle-to infrastructure communication limitations, e.g., bandwidth, dropouts, errors and/or delays.} Ongoing work includes further exposition into the the existence of the optimal solution under different constraint combinations, and the consideration of the terminal speed constrained formulation. {Future work should also address the implementation of the proposed framework under different penetration rates of CAVs and imperfect communication.}
\newpage
\textbf{Appendix}
\begin{appendix}
\section{Proof of Lemma \ref{lem:6}}\label{app:3}
If $v_i(t)-v_{max} \le 0$ becomes active, we have an unconstrained arc (with constant parameters $a_i^{(1)}$,$b_i^{(1)}$,$c_i^{(1)}$,$d_i^{(1)}$) followed by a constrained arc (with constant parameters $a_i^{(2)},~b_i^{(2)},~c_i^{(2)},~d_i^{(2)}$) pieced together at the junction point $t = \tau_s^*$. 
The constrained arc yields at $t = \tau_s^*$ and $t = t_i^m$, 
\begin{subequations}
\begin{align}
   & a_i^{(2)} \cdot \tau_s^* + b_i^{(2)} = 0, \label{eq:tau_s_4}\\
   & a_i^{(2)} \cdot t_i^m + b_i^{(2)} = 0, \label{eq:tau_s_5}\\
   & \frac{1}{2} a_i^{(2)} \cdot (t_i^m)^2 + b_i^{(2)}\cdot (t_i^m) + c_i^{(2)} = v_{max},\label{eq:tau_s_6}\\
   & \frac{1}{6}a_i^{(2)} \cdot (\tau_s^*)^3 + \frac{1}{2}b_i^{(2)}\cdot(\tau_s^*)^2 +  c_i^{(2)}\cdot (\tau_s^*) + d_i^{(2)} \nonumber \\
   & + v_{max}\cdot(t_i^m - \tau_s^*) = p_i(t_i^m). \label{eq:tau_s_7}
\end{align}
\end{subequations}
From \eqref{eq:tau_s_4} and \eqref{eq:tau_s_5}, we have $a_i^{(2)} = 0$ and $b_i^{(2)}=0$. Substituting in \eqref{eq:tau_s_6}, we have $c_i^{(2)} = v_{max}$. Finally, from \eqref{eq:tau_s_7} we have $~d_i^{(2)}=(p_i(t_i^m) - v_{max}\cdot t_i^m)$. The unconstrained arc at the initial condition $t = t_i^0$ yields the following equations: {$\frac{1}{2} a_i^{(1)} \cdot (t_i^0)^2 + b_i^{(1)}\cdot (t_i^0) + c_i^{(2)} = v_i(t_i^0)\label{eq:tau_s_8}, ~
        \frac{1}{6}a_i^{(1)} \cdot (t_i^0)^3 + \frac{1}{2}b_i^{(1)}\cdot(t_i^0)^2 + c_i^{(1)}\cdot (t_i^0) + d_i^{(1)} = p_i(t_i^0).\label{eq:tau_s_9}$
Solving the above two equations by considering $ t_i^0 = 0$, without loss of generality, we have $c_i^{(1)} = v_i(t_i^0)$ and $d_i^{(1)}=0$.} At $\tau_s^*$, we have the following set of equations for the unconstrained arc,
\begin{subequations}
\begin{align}
&    a_i^{(1)} \cdot \tau_s^* + b_i^{(1)} = 0, \label{eq:tau_s_1}\\
 &   \frac{1}{2} a_i^{(1)} \cdot (\tau_s^*)^2 + b_i^{(1)} \cdot \tau_s^* + (v_i(t_i^0) - v_{max}) = 0, \label{eq:tau_s_2}\\
 &   \frac{1}{6} a_i^{(1)} \cdot (\tau_s^*)^3 + \frac{1}{2} b_i^{(1)} \cdot (\tau_s^*)^2 \nonumber\\  
  &  + (v_i(t_i^0) - v_{max})\cdot \tau_s^* - (p_i(t_i^m) - v_{max}\cdot t_i^m)=0.\label{eq:tau_s_3}
\end{align}
\end{subequations}
Substituting $\tau_s^* = -\frac{b_i^{(1)}}{a_i^{(1)}}$ from \eqref{eq:tau_s_1} in \eqref{eq:tau_s_2}, we have $\frac{(b_i^{(1)})^2}{a_i^{(1)}} = 2(v_i(t_i^0)-v_{max})$.
Substituting $\tau_s^* = -\frac{b_i^{(1)}}{a_i^{(1)}}$ from \eqref{eq:tau_s_1} in \eqref{eq:tau_s_3}, we have $    \frac{1}{3} \frac{(b_i^{(1)})^3}{(a_i^{(1)})^2} + \frac{(b_i^{(1)})}{a_i^{(1)}} \cdot (v_{max}-v_i(t_i^0)) - (p_i(t_i^m) - v_{max}\cdot t_i^m) = 0$.
From the last two equations, we obtain $\tau_s^* = - \frac{3 (p_i(t_i^m) - v_{max} \cdot t_i^m)}{(v_{max} - v_i(t_i^0))}$,
where $\tau_s^*$ is an explicit function of the known parameters $p_i(t_i^m), ~v_{max}, ~ v_i(t_i^0) \text{ and } t_i^m$.

\section{Proof of Lemma \ref{lem:7}}\label{app:4}
If $u_i(t)-u_{max} \le 0$ becomes active, we have a constrained arc (with constant parameters $a_i^{(1)},~b_i^{(1)},~c_i^{(1)},~d_i^{(1)}$) followed by an unconstrained arc (with constant parameters $a_i^{(2)},~b_i^{(2)},~c_i^{(2)},~d_i^{(2)}$) pieced together at the junction point $t = \tau_c^*$. Solving \eqref{eq:case2_110} and \eqref{eq:case2_112}-\eqref{eq:case2_114}, we have $    a_i^{(2)} = - \sqrt{\frac{(u_{max})^3}{3(t_i^m)^2\cdot u_{max} + 6t_i^m\cdot v_i(t_i^0)-6L}} \label{eq:lemma7_1}$, where $L= p_i(t_i^m)- p_i(t_i^0)$.
From \eqref{eq:case2_110} and \eqref{eq:case2_111}, $    \tau_c^* = \frac{u_{max}}{a_i^{(2)}}+t_i^m. \label{eq:lemma7_2}$
Finally, substituting $a_i^{(2)}$ into the last equation, the junction point $\tau_c^*$ is given by $     \tau_c^* = t_i^m- \frac{u_{max}}{\sqrt{\frac{(u_{max})^3}{3(t_i^m)^2\cdot u_{max} + 6t_i^m\cdot v_i(t_i^0)-6L}}},$ and can be simplified to  $ \tau_c^* = t_i^m-  \sqrt{\frac{3(t_i^m)^2\cdot u_{max} + 6t_i^m\cdot v_i(t_i^0)-6L}{u_{max}}}$, which is an explicit function of the known boundary parameters $t_i^m,~ p_i(t_i^m),~ p_i(t_i^0),~ v_i(t_i^0), \text{ and } u_{max}.$

\section{Proof of Lemma \ref{lem:8}}\label{app:5}
If $u_i(t)-u_{max} \le 0$ becomes active, we have a constrained arc with constants of integration $a_i^{(1)},~b_i^{(1)},~c_i^{(1)},~d_i^{(1)}$ followed by an unconstrained arc with constants of integration $a_i^{(2)},b_i^{(2)},c_i^{(2)},d_i^{(2)}$, pieced together at the junction point $t = \tau_c^*$. If $v_i(t)-v_{max} \le 0$ becomes active, we have a constrained arc with constants of integration $a_i^{(2)},b_i^{(2)},c_i^{(2)},d_i^{(2)}$ followed by a constrained arc with constants of integration $a_i^{(3)},~b_i^{(3)},~c_i^{(3)},~d_i^{(3)}$ pieced together at the junction point $t = \tau_s^*$. Solving  \eqref{eq:case2_sol1}-\eqref{eq:case2_sol4} for the control constrained arc with $t_i^0=0$, we have $a_i^{(1)} = 0$, $b_i^{(1)} = u_{max}$, $c_i^{(1)} = v_i(t_i^0)$ and $d_i^{(1)} = p_i(t_i^0)$. Solving \eqref{eq:case3_124}-\eqref{eq:case3_123} for the state constrained arc, considering $t_i^0 = 0$ without loss of generality, we have $a_i^{(3)} = 0$, $b_i^{(3)} = 0$, $c_i^{(3)} = v_{max}$  and $d_i^{(3)} = p_i(t_i^m) - v_{max}\cdot t_i^m$. From \eqref{eq:case3_119} and \eqref{eq:case3_122}, we have $\tau_c^* = \frac{u_{max}-b_i^{(2)}}{a_i^{(2)}}$ and $\tau_s^* = -\frac{b_i^{(2)}}{a_i^{(2)}}$ respectively. Substituting the latter into \eqref{eq:case3_120}, \eqref{eq:case3_121}, \eqref{eq:case3_126} and \eqref{eq:case3_128}, and solving the system of equations,{ we have $    a_i^{(2)} = -u_{max}^2\cdot \sqrt{-\frac{1}{\phi}}, \label{eq:lemma8_1} \text{ and }
    b_i^{(2)} = \frac{u_{max}(-2v_i(t_i^0)\sqrt{-\frac{1}{\phi}}+2v_{max}\sqrt{-\frac{1}{\phi}}+1)}{2}, \label{eq:lemma8_2}$
where, $\phi(t_i^m, p_i(t_i^m), v_i(t_i^0), u_{max}, v_{max}) = -24(t_i^m\cdot u_{max}\cdot v_{max} - p_i(t_i^m)\cdot u_{max} + v_i(t_i^0) \cdot v_{max})+12(v_i^2(t_i^0) + v_{max}^2)$.}
Substituting the last results into \eqref{eq:case3_119} and \eqref{eq:case3_122}, the junction points $\tau_s^*$ and $\tau_c^*$ are given as explicit functions of the known parameters $t_i^m,~ p_i(t_i^m),~ v_i(t_i^0),~ u_{max}$ $\text{ and } v_{max}$.

\section{Proof of Lemma \ref{lem:9}}\label{app:1}
From \eqref{eq:jump_case4_3}, we have 
\begin{gather}
 \frac{1}{2}u_i^2(\tau_s^-) + \lambda_i^p(\tau_s^-)\cdot v_i(\tau_s^-)
    + \lambda_i^v(\tau_s^-)\cdot u_i(\tau_s^-)\nonumber\\
    +  \eta_i^d(\tau_s^-)\cdot  (v_{min}-v_i(\tau_s^-)) =     \frac{1}{2}u_i^2(\tau_s^+) + \lambda_i^p(\tau_s^+)\cdot v_i(\tau_s^+)\nonumber\\
    + \lambda_i^v(\tau_s^+)\cdot u_i(\tau_s^+)+ \eta_i^d(\tau_s^+)\cdot (v_{min}-v_i(\tau_s^+)).\label{eq:lem9}
    \end{gather}
Since $v_i(\tau_s^+)=v_{min}$ and $u_i(\tau_s^+)=0$, and from the continuity of state \eqref{eq:speed-continuity} and $\lambda_i^p$ \eqref{eq:jump_case4_1}, we have $\lambda_i^p(\tau_s^-) \cdot v_i(\tau_s^-) = \lambda_i^p(\tau_s^+)\cdot v_i(\tau_s^+)$. From \eqref{eq:17c}, we have $\eta_i^d(\tau_s^-)\cdot (v_{min}-v_i(\tau_s^-))= \eta_i^d(\tau_s^+)\cdot (v_{min}-v_i(\tau_s^+)) =0 $. Hence, \eqref{eq:lem9} reduces to $\frac{1}{2}u_i^2(\tau_s^-)
    + \lambda_i^v(\tau_s^-)\cdot u_i(\tau_s^-) =0$, which implies that either $u_i(\tau_s^-) = 0$ or $\frac{1}{2}u_i(\tau_s^-) + \lambda_i^v(\tau_s^-) = 0$, or both. Since the second term can not hold, we have $u_i(\tau_s^-)=u_i(\tau_s^+)=0$.

\section{Proof of Lemma \ref{lem:11}}\label{app:2}
Since $H_i(\tau_c^+)= H_i(\tau_c^-)$, we have
$\frac{1}{2}u_i^2(\tau_c^-) + \lambda_i^p(\tau_c^-)\cdot v_i(\tau_c^-)
    + \lambda_i^v(\tau_c^-)\cdot u_i(\tau_c^-)+
    \mu_i^b(\tau_c^-)\cdot  (u_{min}-u_i(\tau_c^-)) =\frac{1}{2}u_i^2(\tau_c^+) + \lambda_i^p(\tau_s^+)\cdot v_i(\tau_c^+)
    + \lambda_i^v(\tau_c^+)\cdot u_i(\tau_c^+)+ \mu_i^b(\tau_c^+)\cdot (u_{min}-u_i(\tau_c^+))$.
From the continuity of the state  \eqref{eq:speed-continuity} and $\lambda_i^p$ at $t=\tau_c$, we have $\lambda_i^p(\tau_c^-)\cdot v_i(\tau_c^-) = \lambda_i^p(\tau_c^+)\cdot v_i(\tau_c^+)$. From \eqref{eq:kkt1} we have $ \mu_i^b(\tau_c^-) \cdot (u_{min}-u_i(\tau_c^-))=\mu_i^b(\tau_c^+)\cdot (u_{min}-u_i(\tau_c^+))=0$. After simplifying, we have either $u_i(\tau_c^+) = u_i(\tau_c^-)$ or $\frac{1}{2}(u_i(\tau_c^+) + u_i(\tau_c^-)) + \lambda_i^v(\tau_c^+)=0$. Both the equations lead to the continuity in control input $u_i(t)$ at time $t=\tau_c$, i.e., $u_i(\tau_c^+) = u_i(\tau_c^-)$.
\end{appendix}
\bibliographystyle{abbrvnat}
\bibliography{tcst_constraint_ref, IDS_Publications_03212021}
\end{document}